\documentclass{amsart}
\usepackage{graphicx}
\usepackage{amscd}
\usepackage{amsmath}
\usepackage{amsfonts}
\usepackage{amssymb}
\newtheorem{theorem}{Theorem}
\theoremstyle{plain}

\newtheorem{corollary}{Corollary}

\newtheorem{definition}{Definition}
\newtheorem{example}{Example}

\newtheorem{lemma}{Lemma}
\newtheorem{notation}{Notation}

\newtheorem{proposition}{Proposition}
\newtheorem{remark}{Remark}

\numberwithin{equation}{section}

\begin{document}
\title{Thickness Formula and C$^{1}$-Compactness \\for C$^{1,1}$ Riemannian Submanifolds}
\author{Oguz C. Durumeric}
\address{Department of Mathematics\\
University of Iowa\\
Iowa City, Iowa 52242}
\email{odurumer@math.uiowa.edu}
\thanks{}
\date{January 14, 2002, revised March 2002}
\dedicatory{}

\begin{abstract}
The properties of normal injectivity radius $i(K,M)$ (thickness), of $C^{1,1}
$ submanifolds $K$ of complete Riemannian manifolds $M$ are studied. We
introduce the notion of \ geometric focal distance for $C^{1,1}$ submanifolds
by using metric balls. A formula for $i(K,M)$ in terms of the double critical
points and the geometric focal distance is proved. The thickness of knots and
ideal knots relate to the study of DNA\ molecules and other knotted polymers.
We prove that the set of all $C^{1,1}$ submanifolds $K$ of a fixed manifold
$M$ contained in a compact subset $D\subset M$ and $i(K,M)\geq c>0$ is
$C^{1}-$compact and this collection has finitely many diffeomorphism and
isotopy types. Estimates on upper bounds for the number of such types are
constructible, and we calculate them for submanifolds of $\mathbf{R}^{n}.$
$C^{1}-$compactness is related to Gromov's compactness theorem, but it is an
extrinsic and isometric embedding type theorem.
\end{abstract}\maketitle

\section{Introduction}

Let $M^{n}$ denote a complete connected $n$-dimensional Riemannian manifold.
For a compact $k$-dimensional $C^{1}$ submanifold $K^{k}$ $(\partial
K=\emptyset)$ of $M^{n}$, the normal exponential map, $\exp^{N}$ on the normal
bundle of $K$ in $M$ and its normal injectivity radius $i(K,M)$ are well
defined. If $K$ is $C^{1,1}$, then $i(K,M)>0.$ We will introduce the notion of
''Geometric Focal Distance'' by using metric balls, which naturally extends
the notion of the focal distance of smooth category to $C^{1}$ category in
Riemannian manifolds. We prove a formula for $i(K,M)$ in terms of geometric
focal distance and double critical points for $C^{1,1}$ submanifolds, and that
the set of all submanifolds $K$ of a fixed manifold $M $ contained in a
compact subset $D\subset M$ and $i(K,M)$ bounded away from zero is $C^{1}%
-$compact. These results are essential to the study of the maximization of
$i(K,M)$. The motivation for the maximization of $i(K,M)$ comes from two
directions- the ideal knots and the history of maximization of the intrinsic
injectivity radius.

The thickness of a knotted curve is the radius of the largest tubular
neighborhood around the curve without intersections of normal discs, that is
$i(K,M).$ The ideal knots are the embeddings of $S^{1}$ into $\mathbf{R}^{3},
$ maximizing $i(K,M)$ in a fixed isotopy (knot) class of fixed length. As
noted in [Ka], ''...the average shape of knotted polymeric chains in thermal
equilibrium is closely related to the ideal representation of the
corresponding knot type''. ''Knotted DNA molecules placed in certain solutions
follow paths of random closed walks and the ideal trajectories are good
predictors of time averaged properties of knotted polymers'' as a biologist
referee pointed out to the author. The analytical properties ideal knots will
be tools in the research on the physics of knotted polymers. Theorem I and the
methods developed in this article are used extensively in [D6] where we study
the local structure ideal knots in $\mathbf{R}^{3}.$

Studying ideal knots in $\mathbf{R}^{3}$ corresponds to placing molecules in
homogenous solutions with uniform conditions. Studying ideal knots in
Riemannian manifolds, i.e. varying metrics, may bring new possibilities with
varying conditions, such as inhomogeneous solutions.

For a compact Riemannian manifold $M$, let $d(M)$, $v(M)$ and $i(M)$ denote
its (intrinsic) diameter, volume, and injectivity radius of its exponential
map, respectively. Maximization of $i(M)$ for fixed $d(M)$ or $v(M)$ has a
long history. $i(M)\leq d(M)$ and equality holds if and only if $M$ is a
Blaschke manifold, Warner [Wa], Besse [Be]. It is conjectured that a Blaschke
manifold is isometric to a sphere or a projective space with the standard
metrics up to rescaling of the metric. Berger proved that $v(M)/i(M)^{n}\geq
v(S^{n}(1))/i(S^{n}(1))^{n},$ and equality holds if and only if $M^{n}$ is
isometric to a standard sphere $S^{n}(r)$, by using an inequality proved by
Kazdan. This resolves the $S^{n}$ and $\mathbf{RP}^{n}$ cases of the Blaschke
conjecture. See Besse [Be] and Berger [B] for the literature on Blaschke
manifolds as well as proofs by Berger and Kazdan.

As well as finding these ideal metrics, we examine the topological
restrictions imposed by large injectivity radii. The class defined by
$i(M)^{n}/v(M)\geq c_{1}>0$ is Hausdorff-Gromov precompact, see [Gr, Prop
5.2], [GWY] and [Cr, prop 14]. However, the condition $i(M)/d(M)\geq c_{2}>0$
does not provide precompactness without curvature restrictions. By the
author's work [D4], and [Y], one can estimate a priori upper bounds for the
number of possible homotopy types for $M$ and Betti numbers [D4], and
fundamental group [D5] in terms of $c_{1}.$ The author also studied manifolds
with large injectivity radii: $c_{2}\approx1$ in [D1, D2, D3].

In this article, we approach the normal injectivity radius $i(K,M)$ from a
more general point of view. Let $\mathcal{D}^{\infty}(k,\varepsilon,D;M)$
$=\{(K,M):K\in C^{\infty},$ $\dim K=k,$ $K\subset D,$ and $i(K,M)\geq
\varepsilon\}$ where $D$ is a compact subset of $M$. The behavior of focal
points and $i(K,M)$ are better understood in the smooth category, but
$\mathcal{D}^{\infty}$ is not complete under $C^{1}$ topology. Since
$i(K,M)\geq\varepsilon$ restricts curvature in a sense, the completion must
include $C^{1,1}$ submanifolds. By Proposition 2, $i(K,M)$ is an upper
semi-continuous function of $K$ in $C^{1}$ topology: $\lim\sup_{m}$
$i(K_{m},M)\leq i(K_{\infty},M).$ The extremal cases are more likely not to
occur in $\mathcal{D}^{\infty}$. Very few ideal knots in $\mathbf{R}^{3}$ are
expected to be $C^{2}$, and possibly the unknotted standard circles are the
only ones. This requires the study of $i(K,M)$ in $C^{1,1}$ category.

The formal definitions will be given in section 2. $F_{g}(K)$ is the geometric
focal distance defined in terms of local intersections with metric balls,
$MDC(K)$ is the length of the shortest geodesic normal to $K$ at both of its
endpoints on $K,$ and the ''rolling bead/ball radius, $R_{O}(K,M)$'' is the
largest radius of open metric balls which are tangent to $K$ without
intersecting $K$ elsewhere. We prove the following expected formula for
$i(K,M)$.

\begin{theorem}
For every complete connected smooth Riemannian manifold $M$ and every compact
$C^{1,1}$ submanifold $K$ $(\partial K=\emptyset)$ of $M,$
\[
i(K,M)=R_{O}(K,M)=\min\{F_{g}(K),\frac{1}{2}MDC(K)\}.
\]
\end{theorem}

For a $C^{1,1}$ curve $\gamma,$ $\gamma^{\prime\prime}$ and the curvature
$\kappa\gamma$ of $\gamma$ exist almost everywhere by Rademacher's Theorem.
The supremum of $\kappa\gamma$ is taken on the set of all points where
$\kappa\gamma$ exists. See [D6], Lemma 2 for a proof of $F_{k}(\gamma
)=F_{g}(K^{1})$ in $\mathbf{R}^{n}$. We prove the following corollary for any
dimensions $n>k\geq1,$ for $K^{k}\subset\mathbf{R}^{n}$ in Proposition 12.

\begin{corollary}
(Thickness Formula for Curves in $\mathbf{R}^{n}$) For every simple$,$
$C^{1,1}-$closed curve $\gamma$ in $\mathbf{R}^{n}$ and $K=image(\gamma)$,
\[
i(K,M)=R_{O}(K,M)=\min\{F_{k}(\gamma),\frac{1}{2}MDC(K)\},
\]
where $F_{k}(\gamma)=\left(  \sup\kappa\gamma\right)  ^{-1}$.
\end{corollary}

The formula $i(K,M)=\min\{F_{g}(K),MDC(K)/2\}$ was proved for $C^{2}-$knots in
$\mathbf{R}^{3}$ in [LSDR], and for $C^{1,1}$-knots in $\mathbf{R}^{3}$ by
Litherland in [L]. Nabutovsky [N] extensively studied $C^{1,1}$ hypersurfaces
$K^{n-1}$ in $\mathbf{R}^{n}$ and their injectivity radii. Some of our results
overlap with [N] in this special case$.$ [N] proves the upper semicontinuity
of $i(K^{n-1},\mathbf{R}^{n})$, the lower semicontinuity of
$v(K)/i(K,\mathbf{R}^{n})^{n-1}$ in $C^{1}$ topology, and the compactness of
the class of hypersurfaces with $i(K,\mathbf{R}^{n})$ bounded from below$.$
The analogous (codimension 2) results were obtained by Litherland in [L] for
$C^{1,1}$ knots in $\mathbf{R}^{3}$. Their proofs use $\varepsilon
-$approximations or curvature, while ours use intersections with metric balls.
The relations between curvature and $F_{g}(K)$ are simple in all spaces of
constant curvature. The equality $i(K,M)=R_{O}(K,M),$ a rolling ball/bead
description of the injectivity radius in $\mathbf{R}^{n}$, was known by
Nabutowsky for hypersurfaces, and by Buck and Simon for $C^{2}$ curves, [BS].
The notion of the global radius of curvature developed by Gonzales and
Maddocks for smooth curves in $\mathbf{R}^{3}$ defined by using circles
passing through 3 points of the curve in [GM] is a different characterization
of $i(K,\mathbf{R}^{3})$ from $R_{O}$ due to positioning of the circles and
metric balls.

Gromov's Compactness Theorem was first stated in [Gr] and some of its details
were clarified by Katsuda in [K]. The $C^{1,\alpha}$ estimates of the metrics
of bounded curvature in harmonic coordinates by Jost and Karcher [JK] were
used to complete the proof by Peters [P] and Greene and Wu [GW]. [P] and [GW]
proved the optimal a priori $C^{1,\alpha}$ regularity of the limit metric and
obtained the Lipschitz convergence intrinsically by studying the transition
functions. Gromov's proof, the clarifications by Katsuda, and the version by
Pugh [Pu] rely on Whitney type, non-isometric embeddings into $\mathbf{R}^{N}$
(large $N$) to show Lipschitz closeness of manifolds.

Let $\mathcal{D}(k,\varepsilon,D;M)=\{(K,M):K\in C^{1,1},$ $\dim K=k,$ $K$ is
connected, $K\subset D,$ and $i(K,M)\geq\varepsilon\},$ where $(K,M)$ denotes
an embedding $e:K\rightarrow(M,g_{0})$ with the induced submanifold metric
$e^{\ast}g_{0}$ on $K.$

\begin{theorem}
For every complete connected Riemannian manifold $(M^{n},g_{0})$, every
compact subset $D\subset M$, $1\leq k\leq n-1,$ and $\varepsilon>0,$ the
following holds.

i. $\mathcal{D}(k,\varepsilon,D;M)$ has finitely many diffeomorphism and
isotopy classes.

ii. $\mathcal{D}(k,\varepsilon,D;M)$ is sequentially compact in $C^{1}%
$-topology, i.e. every sequence $\{(K_{m},M)\}_{m=1}^{\infty}$ in
$\mathcal{D}(k,\varepsilon,D;M)$ has $C^{1}$-convergent subsequence whose
limit $(K_{0},M)$ is in $\mathcal{D}(k,\varepsilon,D;M).$

iii. For every $(K,M)\in\mathcal{D}(k,\varepsilon,D;M),$ the induced
submanifold metric $e^{\ast}g_{0}$ is a priori $C^{0,1}$. However, there
exists an isometric $C^{1,1}$ embedding of $(K,g_{\infty})$ onto $(K,e^{\ast
}g_{0})$ in $M$ such that $g_{\infty}$ is $C^{1,\alpha}$ $(\alpha<1)$ in
harmonic coordinates of $K$ where $(K,g_{\infty})$ is a limit of $C^{\infty}$
Riemannian metrics of bounded curvature and injectivity radii with respect to
Lipschitz distance and it is a $C^{1,\alpha}$ Alexandrov space with a well
defined exponential map.
\end{theorem}

In \emph{Theorem 1} and \emph{Section 3}, $K$ is not assumed to be connected.
However, \emph{Theorem 2} is proved for connected $K$, and the general case is
discussed as \emph{Corollary 3} of \emph{Section 4.2}. Theorem 2 is an
extrinsic and isometric embedding type Gromov compactness theorem, but it
differs from the versions above in several aspects. Its proof uses the
intrinsic versions [P] and [GW], and the harmonic coordinates of [JK] to
secure the isometric embedding of the intrinsic limit. We will also show that
there exists $\rho_{0}(k,\varepsilon,D,M)>0$ such that for all $K,L\in
\mathcal{D}(k,\varepsilon,D;M)$ satisfying $L\subset B(K,\rho_{0};M)$ there
exists a continuous isotopy between $K$ and $L$ through $C^{1,1}$ embeddings
of $L$. Thus, if the submanifolds are close in the Hausdorff topology in $M$,
then they are isotopic and close in $C^{1}$-topology. By Theorems 1 and 2, and
Proposition 2, every isotopy class must have a thickest-$i(K,M)$ maximizing submanifold.

The method using the embeddings into $\mathbf{R}^{N}$ could only obtain a
priori $C^{0,1}$ regularity of the limit metric, see [Pu], in contrast to
$C^{1,\alpha}$ $(\alpha<1)$ regularity obtained by the intrinsic proofs of [P]
and [GW]. Any $C^{1}$ simple closed curve $\gamma$ of length $2\pi$ in a
Riemannian manifold has a $C^{0}$ metric induced by the embedding, while
$\gamma$ is intrinsically isometric to standard smooth $S^{1}$, and the
regularity is lost in the embedding. Part (iii) of Theorem 2 emphasizes the
recovery the possible loss of regularity coming from the embedding. We know
more about the geometry of $C^{1,\alpha}$ Alexandrov spaces [Ni] than
$C^{0,1}$ metrics which do not even admit an exponential map a priori. Note
that not all $C^{0,1}$ Riemannian metrics can be $C^{1,1}$ embeddable into
some smooth $M$ with positive thickness.

Section 3 contains the proof of Theorem 1. For a $C^{1,1}$ submanifold $K,$
the normal exponential map $\exp^{N}$ of $K$ is of class $C^{0,1},$ a priori
differentiable almost everywhere. Hence, the Inverse Function Theorem can not
be used to obtain local diffeomorphisms around regular points, and the
property that the focal points being the singular points of $\exp^{N}$ fails.
In general, the set of focal points may not be closed, and $F_{g}$ is not
semi-continuous, see Example 1. We prove a lower-semicontinuity of the normal
cut value in a certain case in Proposition 7 which is sufficient for Theorem
1. Our main tool is the distance functions from the submanifolds. Despite the
similarities of the main theme to the smooth case, our proof contains many
technical details which are not derivable from the classical lemmas of the
smooth cases.

The compactness is discussed in Section 4. The technical details differ from
the previous section. The thickness controls the directional derivatives
$\left\|  f_{%
\operatorname{u}%
\operatorname{u}%
}\right\|  $ a priori for a graph of a function $f$ locally representing $K.$
For smooth $f$, the Hessian is symmetric and one can control $\left\|
f_{uv}\right\|  $ by the polarization identities. For a $C^{1,1}$ function
$f$, $f_{uv}$ are defined a.e., and $f_{uv}=f_{vu}$ a.e. $\left\|
f_{uv}\right\|  $ are not necessarily uniformly bounded in terms of $\left\|
f_{%
\operatorname{u}%
\operatorname{u}%
}\right\|  $ at a point, see Example 1. To apply Arzela-Ascoli Theorem to a
family of such graphs requires equicontinuity of $f_{u},$ for which one may
wish to use uniform boundedness of $f_{uv}$ on the family. Hence, using
mollifiers is a good way to proceed, and one can obtain that smooth
$\mathcal{D}^{\infty}(k,\varepsilon,D;\mathbf{R}^{n})$ is ''almost'' dense in
$\mathcal{D}(k,\varepsilon,D;\mathbf{R}^{n}).$ For a Riemannian manifold $M$,
$F_{g}(p,K)$ depends on the behavior of the metric of $M$ in the given normal
direction as well as the normal curvature. If \emph{codimension(}$K)\geq2$, it
is possible to have high normal curvatures (if defined) in low ambient
curvature directions and low normal curvatures in high ambient curvature
directions at a given point achieving high $F_{g}(p,K).$ For $C^{1,1}$ $K$ in
general, the ''second fundamental form'' $II_{w}^{K}(v)$ is not defined
everywhere, not a quadratic form (Example 1), and not continuous in $v$ even
at a point. Normal curvatures do not satisfy Euler's formula. An averaging
procedure in small neighborhoods with a $C^{1}$ convergence may not be able to
control the change of normal curvatures in different directions, especially if
the limit is discontinuous. We were able to show the existence of
$\delta(\varepsilon)>0$ for which $\mathcal{D}(k,\varepsilon,D;M)$ is in the
closure of $\mathcal{D}^{\infty}(k,\delta,D_{\varepsilon};M)$.

The remaining parts of the proof of Theorem 2 are straightforward. We show
that the collection of the submanifolds in $\mathcal{D}^{\infty}%
(k,\delta,D;M)$ as a collection of manifolds satisfy the conditions of
Gromov's Compactness Theorem. By following Peter's version [P], every
subsequence has a convergent subsequence whose limit is an (intrinsic)
$C^{1,\alpha}$ Riemannian manifold. We use the harmonic coordinates of [JK] to
apply Arzela-Ascoli Theorem and positive thickness $\delta(\varepsilon)$ to
secure the isometric embedding of the (intrinsic) limit into $M.$ All
constants introduced are constructible in terms of $n,k,\varepsilon,D,$ and
$M.$ In the last section, we calculate some estimates for upper bounds of
isotopy and diffeomorphism types for submanifolds of $\mathbf{R}^{n}$ with
thickness bounded away from $0.$

\section{Basic Definitions}

In this section $M^{n}$ always denotes a complete and smooth Riemannian
manifold, and $K^{k}$ denotes a $C^{1}$ submanifold of $M^{n}.$ We refer to
[CE], [GKM] and [DoC] for basic Riemannian geometry. $TK$, $UTK$, $NK$ and
$UNK$ denote the tangent, unit tangent, normal and unit normal bundle to $K$
in $M$. $\exp^{N}:NK\rightarrow M$ denotes the normal exponential map.

\begin{definition}
i. For any metric space $X$ with a distance function $d$, $B(p,r)=\{x\in
X:d(x,p)<r\}$ and $\bar{B}(p,r)=\{x\in X:d(x,p)\leq r\}.$ For $A\subset X$ and
$x\in X$ define $d(x,A)=\inf\{d(x,a):a\in A\}$ and $B(A,r)=\{x\in
X:d(x,A)<r\}.$ The diameter $d(X)$ of $X$ is defined to be $\sup
\{d(x,y):x,y\in X\}.$ If there is ambiguity, we will use $d_{X}$ and $B(p,r;X).$

ii. For $A\subset M^{n}$ and any curve $\gamma$ in $A,$ the length
$\ell(\gamma)$ is defined with respect to the metric space structure of
$M^{n}.$ For any one-to-one curve $\gamma,\ell_{ab}(\gamma)$ and $\ell
_{pq}(\gamma)$ both denote the length of $\gamma$\ between $\gamma(a)=p$ and
$\gamma(b)=q.$

iii. $v(M)$ denotes the volume of a $C^{1}$ Riemannian manifold $M.$
\end{definition}

\begin{definition}
Let $K$ be a $C^{1}$ submanifold of $M$.

i. Define the thickness of $K$ in $M$ or the normal injectivity radius of
$\exp^{N}$ to be
\[
i(K,M)=\sup(\left\{  0\right\}  \cup\{r>0:\exp^{N}:\{v\in NK:\left\|
v\right\|  <r\}\rightarrow M\text{ is one-to-one}\}).
\]
ii. For any $w\in UNK_{p},$ define the normal cut value in the direction $w$
with respect to $K$ to be $r_{w}=\sup\{r:d(\exp^{N}rw,K)=r\}.$
\end{definition}

\begin{definition}
For a smooth and complete Riemannian manifold $M$ and $p\in M,$ define
pointwise injectivity radius
\[
i(p,M)=\sup\{r>0:\exp:\{v\in TK_{p}:\left\|  v\right\|  <r\}\rightarrow
M\text{ is one-to-one}\}.
\]
For a compact subset $D$ of $M,$ define $i(D)=\underset{p\in D}{\min\text{
\ }}i(p,M)>0.$
\end{definition}

\begin{definition}
Let $K$ be a $C^{1}$ submanifold of $M.$ $K$ and $M$ will be suppressed, if
there is no ambiguity. For any $v\in UNK_{p}$ and any $r,$ define

i. $O_{p}(v,r;M)=%
{\displaystyle\bigcup\limits_{w\in v^{\bot}(1)}}
B(\exp_{p}rw,r)$, where $v^{\bot}(1)=\{w\in UTM_{p}:\left\langle
v,w\right\rangle =0\}$

ii. $O_{p}(r;K)=%
{\displaystyle\bigcap\limits_{v\in UTK_{p}}}
O(v,r;K)=%
{\displaystyle\bigcup\limits_{w\in UNK_{p}}}
B(\exp_{p}rw,r)$

iii. $O(r;K)=%
{\displaystyle\bigcup\limits_{p\in K}}
O_{p}(r;K)$

iv. $O_{p}^{c}(v,r;K)=\mathbf{M}-O_{p}(v,r;K).$
\end{definition}

\begin{definition}
Let $K$ be a $C^{1}$ submanifold of $M.$ Define

i. The ball radius of $K$ in $M$ to be $R_{O}(K,M)=\inf\{r>0:O(r;K)\cap
K\neq\emptyset\}$

ii. The pointwise geometric focal distance $F_{g}(p)=\inf\{r>0:p\in
\overline{O_{p}(r;K)\cap K}\}$ for any $p\in K$ and the geometric focal
distance $F_{g}(K)=\inf_{p\in K}F_{g}(p).$
\end{definition}

\begin{definition}
Let $K$ be a $C^{1}$ submanifold of $M.$ A pair of points $p$ and $q$ in $K$
are called a double critical pair for $K,$ if there is a geodesic $\gamma
_{pq}$ of positive length from $p$ to $q$, normal to $K$ at both $p$ and $q,$
and minimal up to its midpoint from both $p$ and $q.$ Define the minimal
double critical distance

$MDC(K)=\inf\{\ell(\gamma_{pq}):\{p,q\}$ is a double critical pair for $K\}.$
\end{definition}

\begin{definition}
By furnishing the Grassmanian bundle $G_{n,k}(TM)$ with a fixed Riemannian
metric, for $C^{1}-$diffeomorphic compact k-dimensional submanifolds $K$ and
$L$ of $M,$ one defines $d_{C^{1}}(K,L)$ to be

$\inf\{\sup_{x\in K}\left(  d_{M}(x,\psi(x))+d_{G}(TK_{x},TL_{\psi
(x)})\right)  :\forall C^{1}-$diffeomorphisms $\psi:K\rightarrow L\}.$
\end{definition}

\section{Thickness Formula}

Throughout this section, we assume that $K$\ is a compact $C^{1,1}%
$\ submanifold of a complete connected smooth Riemannian manifold $M$ and
$\partial K=\emptyset$, unless stated otherwise. $K$ is not assumed to be connected.

\begin{proposition}
$i(K,M)=R_{O}(K,M).$
\end{proposition}

\begin{proof}
i. Choose any $R>i(K,M).$ There exists $p_{i}\in K$ and $v_{i}\in NK_{p_{i}}$
for $i=1,2$ such that $v_{1}\neq v_{2},$ $q=\exp_{p_{1}}^{N}v_{1}=\exp_{p_{2}%
}^{N}v_{2},$ and $\left\|  v_{1}\right\|  \leq\left\|  v_{2}\right\|  <R, $ by
definition of $i(K,M).$ Let $q_{2}=\exp_{p_{2}}^{N}v_{2}R/\left\|
v_{2}\right\|  .$ Since $\gamma_{1}(t)=\exp_{p_{1}}^{N}tv_{1}$ and $\gamma
_{2}(t)=\exp_{p_{2}}^{N}tv_{2}$ are distinct geodesics, $\pi>\measuredangle
(-\gamma_{1}^{\prime}(q),-\gamma_{2}^{\prime}(q)).$ By First Variation, [CE,
GKM], $d(q_{2},p_{1})<d(q_{2},q)+d(q,p_{1})\leq R-\left\|  v_{2}\right\|
+\left\|  v_{1}\right\|  \leq R.$ Thus, $p_{1}\in B(q_{2},R)\subset O_{q_{2}%
}(R;K)$, $O(R)\cap K\neq\emptyset$ and this concludes $R\geq R_{O}(K,M).$ We
have shown that $\forall R>i(K,M),$ $R\geq R_{O}(K,M).$ This proves that
$i(K,M)\geq R_{O}(K,M).$

ii. Choose any $R>R_{O}(K,M)$ so that $O(R)\cap K\neq\emptyset$. There exists
$p\in K$, $v\in UNK_{p},$ $q=\exp_{p}^{N}vR$ such that $B(q,R)\cap
K\neq\emptyset$. Let $p_{1}\in B(q,R)\cap K$ be any point and $d(q,p_{1}%
)=R-\delta,$ for some $\delta>0.$ Let $q_{1}=\exp_{p}^{N}v(R-\frac{\delta}%
{3}),$ and $p_{2}$ be any closest point of $K$ to $q_{1}.$
\[
R^{\prime}=d(q_{1},p_{2})\leq d(q_{1},p_{1})\leq d(q_{1},q)+d(q,p_{1}%
)\leq\frac{\delta}{3}+R-\delta=R-\frac{2\delta}{3}%
\]
For any normal minimal geodesic $\gamma$ from $p_{2}$ to $q_{1}$,
$w=\gamma^{\prime}(p_{2})\in UNK_{p_{2}}.$ So, $q_{1}=\exp_{p}^{N}%
v(R-\frac{\delta}{3})=\exp_{p_{2}}^{N}wR^{\prime},$ but $v(R-\frac{\delta}%
{3})\neq wR^{\prime}$ since $\left\|  w\right\|  =\left\|  v\right\|  =1.$
Consequently, $\exp^{N}|\{v\in NK:\left\|  v\right\|  <R\}$ is not injective
and $R\geq i(K,M).$ We have shown that $\forall R>R_{O}(K,M),$ $R\geq i(K,M).$
This proves that $i(K,M)\leq R_{O}(K,M).$
\end{proof}

\begin{proposition}
Let $K,K_{j},j\in\mathbf{N,}$ be compact $C^{1}$ submanifolds of a complete
Riemannian manifold $M,$ such that $K_{j}\rightarrow K$ in $C^{1}$ sense. Then
$\underset{j\rightarrow\infty}{\lim\sup}R_{O}(K_{j},M)\leq R_{O}(K,M).$
\end{proposition}

\begin{proof}
Choose arbitrary $R>R_{O}(K,M),$ that is $O(R;K)\cap K\neq\emptyset.$%
\begin{align*}
\exists p  &  \in K\text{ }\exists v\in UNK_{p}\text{ }\exists\varepsilon
>0\text{ }\exists q\in B(\exp_{p}^{N}Rv,R-\varepsilon)\cap K\neq\emptyset\\
\forall j\text{ }\exists p_{j}  &  \in K_{j}\text{ }\exists v_{j}\in
UN(K_{j})_{p_{j}}\text{ such that }(p_{j},v_{j})\rightarrow(p,v),\text{ as
}j\rightarrow\infty\\
\exists j_{0}\forall j  &  \geq j_{0}\text{, }d(\exp_{p_{j}}^{N}Rv_{j}%
,\exp_{p}^{N}Rv)<\frac{\varepsilon}{2}\text{ and }\exists q_{j}\in K_{j}\text{
with }d(q,q_{j})<\frac{\varepsilon}{2}\\
\forall j  &  \geq j_{0},\text{ }d(\exp_{p_{j}}^{N}Rv_{j},q_{j})\leq
d(\exp_{p_{j}}^{N}Rv_{j},\exp_{p}^{N}Rv)+d(\exp_{p}^{N}Rv,q)+d(q,q_{j})<R\\
\forall j  &  \geq j_{0},\text{ }O(R;K_{j})\cap K_{j}\neq\emptyset,\text{that
is: }R_{O}(K_{j},M)\leq R\\
R  &  \geq\underset{j\rightarrow\infty}{\lim\sup}R_{O}(K_{j},M)
\end{align*}

We have shown that if $R>R_{O}(K,M)$ then $R\geq\underset{j\rightarrow\infty
}{\lim\sup}R_{O}(K_{j},M).$ Hence, $R_{O}(K,M)\geq\underset{j\rightarrow
\infty}{\lim\sup}R_{O}(K_{j},M).$
\end{proof}

\begin{proposition}
Let $K_{j},j\in\mathbf{N,}$ be a sequence of $C^{1}$ k-dimensional
submanifolds of a complete Riemannian manifold $M^{n}$, such that
$K_{j}\rightarrow K$ in $C^{1}$ sense, where $K$ is compact. If $\lim\inf
_{j}MDC(K_{j})>0$ then $\lim\inf_{j}MDC(K_{j})\geq MDC(K).$
\end{proposition}

\begin{proof}
We will use the same indices for subsequences. Let $a=\lim\inf_{j}MDC(K_{j}),$
and choose a subsequence with $a=\lim_{j}MDC(K_{j})$ and $\forall
j,MDC(K_{j})>0.$ By compactness of $K_{j}$ and positivity of $MDC(K_{j})$,
there exists a minimal double critical pair $\{p_{j},q_{j}\}$ for $K_{j}$,
$\ell(\gamma_{p_{j}q_{j}})=MDC(K_{j}).$ Since $K$ is compact and $a>0,$ there
exist subsequences $p_{j}\rightarrow p_{0}\in K$, $q_{j}\rightarrow q_{0}\in
K,$ and $\gamma_{p_{j}q_{j}}\rightarrow\gamma_{p_{0}q_{0}}$ in $C^{1}$ sense.
Geodesics converge to geodesics, and normality to submanifolds is preserved
under $C^{1}$ limits. $\{p_{j},q_{j}\} $ is a double critical pair for $K.$%
\[
MDC(K)\leq\ell(\gamma_{p_{0}q_{0}})=\lim_{j}\ell(\gamma_{p_{j}q_{j}})=\lim
_{j}MDC(K_{j})=a
\]
\end{proof}

\begin{proposition}
$R_{O}(K,M)\leq\min\{F_{g}(K),\frac{1}{2}MDC(K)\}.$
\end{proposition}

\begin{proof}
This is am immediate consequence of definitions of $F_{g}(K)$ and $MDC(K).$
\end{proof}

\begin{proposition}
Let $v\in UNK_{p}$ be such that $0\leq r_{v}<F_{g}(K)$. Then there are
finitely many and at least two minimal geodesics between $q=\exp_{p}^{N}%
r_{v}v$ and $K.$ Hence, $r_{v}>0.$
\end{proposition}

\begin{proof}
Any geodesic that is a shortest curve between a point of $M-K$ and $K$ is
normal to $K.$ We assume that all geodesics are unit speed and start at $K$
when $s=0.$ Let $\gamma_{0}(s)=\exp_{p}^{N}sv.$ $\forall j\in\mathbf{N}^{+},$
there exists a minimal geodesic $\gamma_{j}$ between $q_{j}=\exp_{p}^{N}%
(r_{v}+\frac{1}{j})v$ and $K.$ Since $\gamma_{j}$ is not minimal between $p$
and $q_{j}$, $\gamma_{0}\neq\gamma_{j},\forall j$. By compactness, and taking
a subsequence and using the same subindices, we can assume that $\gamma
_{j}\rightarrow\gamma_{\infty},$ a minimal geodesic between $q$ and $K.$

Suppose that $\gamma_{0}=\gamma_{\infty}$ or $r_{v}=0.$ Choose $j_{0}$
sufficiently large such that $a=r_{v}+\frac{1}{j_{0}}<F_{g}(K).$ $\forall
j>j_{0},$ $d(q_{j_{0}},\gamma_{j}(0))<d(q_{j_{0}},q_{j})+d(q_{j},\gamma
_{j}(0))\leq a,$ since $\gamma_{j}^{\prime}(q_{j})\neq\gamma_{0}^{\prime
}(q_{j})$ and First Variation. $\gamma_{j}(0)\in B_{a}(q_{j_{0}})\cap K.$
$\gamma_{j}(0)\neq p,\forall j,$ but $\gamma_{j}(0)\rightarrow p$. Hence,
$p\in\overline{K\cap O_{p}(a)}$ which contradicts $a<F_{g}(K).$ This shows
that $r_{v}>0$ and $\gamma_{0}\neq\gamma_{\infty},$ that is there are at least
two geodesics between $q$ and $K.$

Suppose that there are infinitely many minimal geodesics $\theta_{j}$ between
$K$ and $q.$ By compactness, there exists a convergent subsequence of distinct
geodesics $\theta_{j}\rightarrow\theta_{0}$ which is also minimal between $K$
and $q.$ Then one uses a proof similar to above, with $\theta_{j}%
(0)\rightarrow\theta_{0}(0)=p^{\prime},$ to show $p^{\prime}\in\overline{K\cap
O_{p^{\prime}}(a^{\prime})}$ where $a^{\prime}$ is chosen similar to above
$a^{\prime}=r_{\theta_{0}^{\prime}(0)}+\frac{1}{j_{0}}<F_{g}(K).$ Hence, there
are finitely such geodesics.
\end{proof}

\begin{lemma}
Let $p\in K$ be such that $F_{g}(p)>0.$ $\forall v\in UNK_{p},\forall
r<F_{g}(p),$ $q=\exp_{p}^{N}rv,$ there exists an open disc $D$ of $K,$ such
that $p\in D$ and $\forall x\in D-\{p\},$ $d(x,q)>r.$
\end{lemma}

\begin{proof}
Choose $a$ such that $r<a<F_{g}(p).$ Let $\gamma(s)=\exp_{p}^{N}sv,$ and
$q=\gamma(a).$ $p\notin$ $\overline{K\cap O_{p}(a)},$ since $a<F_{g}(p).$
Hence, there exists an open disc $D$ of $K$ such that $p\in D$ and $D\cap
O_{p}(a)=\emptyset.$ $B(q,r)\subset B(q^{\prime},a)$ and $\forall x\in D,$
$d(x,q)\geq r$ and $d(x,q^{\prime})\geq a.$ Let $x\in D$ be such that
$d(x,q)=r.$ If $\gamma_{x}$ is any normal minimal geodesic from $x$ to $q$,
distinct from $\gamma$, that is $\gamma^{\prime}(q)\neq\gamma_{x}^{\prime}%
(q)$, then by First Variation, $d(x,q^{\prime})<d(x,q)+d(q,q^{\prime})\leq a.$
This contradicts $D\cap O_{p}(a)=\emptyset.$ Finally, $\gamma$ and $\gamma
_{x}$ must follow the same minimal geodesic and $x=p.$
\end{proof}

\begin{proposition}
Let $v\in UNK_{p}$ be such that $r_{v}<F_{g}(K)$ and there are two distinct
minimal geodesics $\gamma_{1}$ and $\gamma_{2}$ between $q=\exp_{p}^{N}r_{v}v$
and $K.$ Then either $\measuredangle(\gamma_{1}^{\prime}(q),\gamma_{2}%
^{\prime}(q))=\pi$ or $R_{O}(K,M)<r_{v}.$
\end{proposition}

\begin{proof}
Assume that $\measuredangle(\gamma_{1}^{\prime}(q),\gamma_{2}^{\prime}%
(q))<\pi,$ to show that $R_{O}(K,M)<r_{v}.$ Let $\gamma_{j}(0)=p_{j}\in K, $
for $j=1,2$. $\gamma_{j}(r_{v})=q,$ $d(q,p_{j})=d(q,p)=r_{v},$ for $j=1,2$.
Choose $\varepsilon>0$ small enough that

1. $D_{j}=B(p_{j},\varepsilon)\cap K$ is a small open disc in $K$ and
$\overline{D_{j}}$ is compact, for $j=1,2$,

2. $\overline{D_{1}}\cap\overline{D_{2}}=\emptyset,$ and

3. $\forall x\in\overline{D_{1}}\cup\overline{D_{2}}-\{p_{1},p_{2}\},$
$d(x,q)>r_{v}$ by the previous lemma.

Let $\delta=\min\{d(x,q)-r_{v}:x\in\partial D_{1}\cup\partial D_{2}\}>0.$
Choose $w\in UTM_{q}$ such that $\measuredangle(w,\gamma_{j}^{\prime
}(q))>\frac{\pi}{2},$ for $j=1,2$. By the First Variation, $d(p_{j},\exp
_{q}tw)$ decreases strictly, for small $t>0$ and for $j=1,2$. If $q$ is on
$cutlocus(p_{j})$, then one can use Toponogov's Theorem, see [CE], [GKM].
There exists $t_{0}\in(0,\frac{\delta}{3})$ and $q_{0}=\exp_{q}t_{0}w$ such
that $r_{v}-\frac{\delta}{3}<d(p_{j},q_{0})<d(p_{j},q)=r_{v},$ for $j=1,2.$
Let $m_{j}$ be the closest point of $\overline{D_{j}}$ to $q_{0}$, for $j=1,2.$

Suppose that $m_{j}\in\partial D_{j}$, for $j=1$ or $2.$ Then, we obtain a
contradiction as follows:
\begin{align*}
d(m_{j},q_{0})  &  \geq d(m_{j},q)-d(q,q_{0})\geq r_{v}+\delta-\frac{\delta
}{3}=r_{v}+\frac{2\delta}{3}\\
d(m_{j},q_{0})  &  \leq d(p_{j},q_{0})\leq d(p_{j},q)+d(q,q_{0})\leq
r_{v}+\frac{\delta}{3}%
\end{align*}

Hence, $m_{j}$ are interior points of $D_{j}$, for $j=1$ and $2$. The minimal
geodesics from $m_{j}$ to $q_{0}$ are normal to $K$ at $m_{j}.$ $m_{1}\neq
m_{2},$ since $\overline{D_{1}}\cap\overline{D_{2}}=\emptyset.$ Finally,
$\exp_{K}^{N}$ fails to be injective on the closed disc bundle of radius
$\max(d(p_{1},q_{0}),d(p_{2},q_{0}))<r_{v}.$ $R_{O}(K,M)=i(K,M)<r_{v}.$
\end{proof}

\begin{proposition}
If $r_{v}<F_{g}(K),$ then $\underset{w\rightarrow v}{\lim\inf}r_{w}\geq
r_{v}.$ That is, $r_{v}$ is lower semi-continuous in $v$ on $UNK$ when
$r_{v}<F_{g}(K).$
\end{proposition}

\begin{proof}
Suppose not, and choose $v_{j}\rightarrow v$ such that $\underset
{v_{j}\rightarrow v}{\lim}r_{v_{j}}=L<r_{v}<F_{g}(K),$ where $v\in UNK_{p}$
and $v_{j}\in UNK_{p_{j}},\forall j\in\mathbf{N,}$ and $p_{j}\rightarrow p.$
We will obtain a contradiction in both cases below.

\textbf{Case 1}. $L>0.$ By Proposition 6, $\forall j\in\mathbf{N,}$ there
exists $p_{j}^{\prime}\in K,$ $u_{j}\in UNK_{p_{j}^{\prime}},$ $q_{j}\in M$
such that $u_{j}\neq v_{j}$, $r_{u_{j}}=r_{v_{j}}$ and $\exp_{p_{j}}^{N}%
v_{j}r_{v_{j}}=\exp_{p_{j}^{\prime}}^{N}u_{j}r_{u_{j}}:=q_{j}.$ By taking
subsequences and using same indices, we may assume that $p_{j}^{\prime
}\rightarrow p^{\prime}$, $u_{j}\rightarrow u$ and $q_{j}\rightarrow
q=\exp_{p}^{N}Lv=\exp_{p^{\prime}}^{N}Lu.$ The case of $u\neq v$ cannot occur,
since $L<r_{v}.$ Hence, we need to study the case of $p=p^{\prime}$ and $v=u.$
Let $c_{v}=\sup\{t:d(p,\exp_{p}tv)=t\}$ be the cut value of the exponential
map $\exp:TM\rightarrow M$ in the direction of $v$ from $p$. Obviously,
$c_{v}\geq r_{v}>L=d(p,q),$ from the definition of the normal cut value.

\qquad See [CE, p.93, 95] or [DoC, p267-276], for the $C^{\infty}$ Riemannian
manifolds $M,$ to conclude that

i. $q$ is not conjugate to $p$ along the unique minimal geodesic $\exp_{p}tv,$
and $q\notin cutlocus(p)$ and hence,

ii. $p$ is not conjugate to $q$ along the unique minimal geodesic $\exp
_{p}(L-t)v=\exp_{q}tw,$ $p\notin cutlocus(q)$ and $c_{w}>L.$ By [DoC, p276 or
CE p.94], the cut value function $c_{(.)}:UM\rightarrow\lbrack0,\infty] $ is
continuous and the tangential cutlocus is a closed subset of $UM$.

Hence, there exists $\varepsilon>0$ satisfying:

1. $0<\varepsilon<\frac{1}{2}\min(F_{g}(K)-L,L),$ and

2. $p$ is the unique closest point of $K$ to $q$ in $K\cap B(p,2\varepsilon),$
by Lemma 1, and

3. $\forall x\in B(p,\varepsilon),\forall y\in B(q,\varepsilon),$ $x\notin
cutlocus(y)$ and unique minimal geodesics $\gamma_{xy}$ vary continuously on
$B(p,\varepsilon)\times B(q,\varepsilon).$

Let $K^{\prime}=K\cap B(p,\varepsilon),$ and consider $\partial K^{\prime} $
with respect to the topology of $K$. Define $\delta:=\frac{1}{4}\min
\{d_{M}(x,q)-L:x\in\partial K^{\prime}\}.$ Observe that $0<\delta\leq
\frac{\varepsilon}{4},$ since $p$ is the unique closest point of
$\overline{K^{\prime}}$ to $q$ and triangle inequality. Choose\ and fix
sufficiently large $j_{0}$ with $q_{j_{0}}\in B(q,\delta;M),$ $p_{j_{0}%
},p_{j_{0}}^{\prime}\in B(p,\frac{\delta}{2};K)$ and $r_{u_{j_{0}}%
}=r_{v_{j_{0}}}<F_{g}(K).$ There exists a curve $\gamma$ of length $\leq
\delta$ in $K$ between $p_{j_{0}}$ and $p_{j_{0}}^{\prime}.$ Define
$f(x)=d_{M}(x,q_{j_{0}}):K^{\prime}\rightarrow\mathbf{R}$ and set $m=\min
f=d_{M}(p_{j_{0}},q_{j_{0}})=d_{M}(p_{j_{0}}^{\prime},q_{j_{0}}).$ By triangle
inequality, we have:
\begin{align*}
0  &  <L-2\delta\leq m\leq f(x)\leq L+\delta+\varepsilon<F_{g}(K),\forall x\in
K^{\prime}\\
m  &  \leq f(\gamma(t))\leq m+\frac{\delta}{2}\leq L+\frac{5\delta}{2}%
<F_{g}(K)\text{ and }\\
\underset{x\in\partial K^{\prime}}{\min}f(x)  &  \geq\underset{x\in\partial
K^{\prime}}{\min}d(x,q)-\delta=L+3\delta
\end{align*}

$f\in C^{1}$, since $K^{\prime}\cap(\{q_{j_{0}}\}\cup cutlocus(q_{j_{0}%
}))=\emptyset$ and $K$ is $C^{1,1}.$ A point $x\in K^{\prime}$ is a critical
point of $f$ if and only if the minimal geodesic from $q_{j_{0}}$ to $x$ is
normal to $K.$ All of the critical points of $f$ are isolated strict local
minima by $f(x)<F_{g}(K)$ and Lemma 1. For an isolated local strict minimum
point $x_{0}$, $x_{0}\notin\overline{f^{-1}((0,f(x_{0}))}.$ As $b$ increases,
$f^{-1}((0,b))$ will gain new components at each critical point $x_{0}.$ Away
from critical points $f^{-1}(b)$ is a codimension 1 submanifold with a normal
$\nabla f$ pointing away from $f^{-1}((0,b))$. Hence, as $b$ increases, the
number of components of $f^{-1}((0,b))$ will not decrease at regular points$.$
By Milnor [M, p.12], for $m<b<L+3\delta$, $f^{-1}((0,b))$ is a disjoint union
of open sets where each component is away from $\partial K^{\prime}$ and
contains exactly one local minimum. However, $\gamma\subset f^{-1}%
((0,L+\frac{11\delta}{4}))\subset int(K^{\prime})$ and the end points\ of
$\gamma$, $p_{j_{0}}$ and $p_{j_{0}}^{\prime}$ are the absolute minima of $f.$
This gives a contradiction. Consequently,\ the case of $p=p^{\prime}$ and
$v=u$ can't occur either.

\textbf{Case 2}. $L=0.$ Let $\eta>0$ be the infimum of the pointwise
injectivity radius of $\exp_{p}:TM_{p}\rightarrow M$ where $p$ ranges over
$d(K)$ neighborhood of $K$ in $M.$ Let $p_{j}\rightarrow p$, $p_{j}^{\prime
}\rightarrow p^{\prime}$ and $q_{j}\rightarrow q$ be chosen as in $L>0$ case.
$p=p^{\prime}=q$ since $L=0.$ Choose $j$ sufficiently large so that
$max(d(p_{j},q_{j}),d(p_{j}^{\prime},q_{j}))<\frac{\eta}{4}.$ Suppose that
there exists a curve $\gamma$ in $K$ between $p_{j}$ and $p_{j}^{\prime}$ of
length $\leq\frac{\eta}{2}.$ Apply the method of in Case 1 to the critical
points of $C^{1}$-function $f(x)=d(x,q_{j}):K\cap B(q_{j},\eta)\rightarrow
(0,\eta)$, which are strict local minima to obtain a contradiction. Hence, all
curves in $K$ between $p_{j}$ and $p_{j}^{\prime}$ have length $>\frac{\eta
}{2}$, and hence $d_{K}(p_{j},p_{j}^{\prime})\geq\frac{\eta}{2}$ for
sufficiently large $j.$ That is not possible since $p_{j}\rightarrow p$ and
$p_{j}^{\prime}\rightarrow p^{\prime}=p.$

Finally, we obtained contradictions in both cases. We can conclude that there
exists no subsequence $v_{j}\rightarrow v$ such that $\underset{v_{j}%
\rightarrow v}{\lim}r_{v_{j}}=L<r_{v}.$
\end{proof}

\subsection{Proof of Theorem 1}

\begin{proof}
By Propositions 1 and 4: $i(K,M)=R_{O}(K,M)\leq\min\{F_{g}(K),\frac{1}{2}MDC(K)\}.$

Claim.$\underset{v\in UNK}{\inf}r_{v}=i(K,M)$.

Let$\underset{v\in UNK}{\inf}r_{v}=r.$ Choose any $\rho>r$ and let
$x=\frac{\rho+r}{2}$. There exists $v_{0}\in UNK_{p_{0}}$ such that
$d(\exp_{p_{0}}^{N}xv_{0},K)<x$, hence $\exp_{p_{0}}^{N}tv_{0}$ and any
minimal geodesic $(\neq\exp_{p_{0}}^{N}tv_{0})$ between $\exp_{p_{0}}%
^{N}xv_{0}$ and $K\ $will have a common point, $\exp_{p_{0}}^{N}xv_{0}.$
Hence, injectivity of $\exp^{N}$ fails in $\rho-$neighborhood and
$\rho>i(K,M).$ Thus, $r\geq i(K,M).$

Choose any $\rho>i(K,M)$ and let $x=\frac{\rho+i(K,M)}{2}$. $\exp^{N}$ fails
to be injective in the $x-$neighborhood. $\exp_{p_{1}}^{N}t_{1}v_{1}%
=\exp_{p_{2}}^{N}t_{2}v_{2}$ for $t_{j}<x$ and $v_{1}\neq v_{2}.$ Then
$\exp_{p_{2}}^{N}tv_{1}$ is not minimal to $K$ for $t>t_{1},$ by First
Variation. Hence, $r_{v_{1}}\leq t_{1}<x<\rho.$ Thus $r<\rho$, to conclude
$r\leq i(K,M).$ This proves the claim.

If $i(K,M)=F_{g}(K),$ then there is nothing to prove. Hence, assume that
$\underset{v\in UNK}{\inf}r_{v}=i(K,M)<F_{g}(K).$ Choose $v_{j}\in UNK_{p_{j}%
},\forall j\in\mathbf{N,}$ and $\underset{j}{\lim}r_{v_{j}}=\underset{v\in
UNK}{\inf}r_{v}.$ By compactness of $K,$ there exists $v_{\infty}\in
UNK_{p_{\infty}}$ and a subsequence which we will denote with the same indices
such that $v_{j}\rightarrow v_{\infty}$ and $p_{j}\rightarrow p_{\infty}.$ By
Proposition 7, $\underset{j}{\lim}r_{v_{j}}\geq r_{v_{\infty}}\geq
\underset{v\in UNK}{\inf}r_{v}.$ Hence $R_{O}(K,M)=i(K,M)=\underset{v\in
UNK}{\inf}r_{v}=r_{v_{\infty}}<F_{g}(K).$ By Propositions 5 and 6, there are
two distinct minimal geodesics $\gamma_{1}$ and $\gamma_{2}$ between
$q=\exp_{p_{\infty}}^{N}r_{v_{\infty}}v_{\infty}$ and $K\ $with
$\measuredangle(\gamma_{1}^{\prime}(q),\gamma_{2}^{\prime}(q))=\pi$. This
means that $\frac{1}{2}MDC(K)=i(K,M)$.
\end{proof}

\begin{example}
Let $h(\theta):\mathbf{R}\rightarrow\mathbf{R}$ be smooth with $h(\theta
+\pi)=h(\theta),\forall\theta.$ Consider the graph $K$ of $z=f(x,y)=\frac
{1}{2}r^{2}h(\theta)$ in $\mathbf{R}^{3}$ where $(r,\theta)$ is the polar
coordinates. $\ f\in C^{1,1}.$ Obviously, $f_{%
\operatorname{u}%
\operatorname{u}%
}(0,0)=h(\theta)$, if $u=(\cos\theta,\sin\theta),$ and $f_{yx}(0,0)=\frac
{1}{2}h^{\prime}(0).$ Away from $(0,0),$ $f$ is smooth and
\[
f_{xx}(x,y)=h(\theta)-\frac{xy}{x^{2}+y^{2}}h^{\prime}(\theta)+\frac{y^{2}%
}{2(x^{2}+y^{2})}h^{\prime\prime}(\theta).
\]

a) Choose $h$ such that $h$ is identically $0$ on an open subset of
$\mathbf{R}$ near $\theta=0,$ $\left\|  h\right\|  _{\infty}=1,$ and
$h^{\prime\prime}(\frac{\pi}{2})\approx n,$ large $n.$ $F_{g}(\mathbf{0}%
,K)=1,$ and every neighborhood of $\mathbf{0}$ contains planar points
$\mathbf{a}$ with $F_{g}(\mathbf{a},K)=\infty,$ as well as points $\mathbf{b}$
with $F_{g}(\mathbf{b},K)\approx\frac{2}{n}.$ Hence, $F_{g}$ is not a
semicontinuous function into $[0,\infty].$ Same is true for the normal
cutvalue if $F_{g}$ is the controlling factor.

b) Choose $h$ such that $\left\|  h\right\|  _{\infty}=1,$ and $h^{\prime
}(0)\approx n,$ large $n,$ to observe that $\forall u,\left\|  f_{%
\operatorname{u}%
\operatorname{u}%
}(0,0)\right\|  \leq1,$ and $F_{g}(\mathbf{0},K)=1$ does not control $\left\|
f_{yx}(0,0)\right\|  .$
\end{example}

\section{Compactness}

Throughout this section we will assume the following. $M$ denotes a smooth
connected complete $n$-dimensional Riemannian manifold, $D\subset M$ denotes a
compact subset, and $K$ denotes a $k$-dimensional compact \textbf{connected}
(except in Section 4.2) $C^{1,1}$ manifold. $(K,M)$ denotes that $K$ is a
Riemannian submanifold with a particular embedding and furnished with the
induced submanifold metric. $(K,g)$ denotes a manifold with a metric $g$
without any indication of any embedding. We refer to DoCarmo [DoC], for basic
submanifold theory for Riemannian manifolds.

\begin{definition}
$\ $

$\mathcal{A}(k,\varepsilon,D;M)=\{(K,M):K\in C^{1,1},$ $\dim K=k,$ $K\subset
D,$ and $i(K,M)>\varepsilon\}$

$\mathcal{A}^{\infty}(k,\varepsilon,D;M)=\{(K,M):K\in C^{\infty},$ $\dim K=k,$
$K\subset D,$ and $i(K,M)>\varepsilon\}$

$\mathcal{D}(k,\varepsilon,D;M)=\{(K,M):K\in C^{1,1},$ $\dim K=k,$ $K\subset
D,$ and $i(K,M)\geq\varepsilon\}$

$\mathcal{D}^{\infty}(k,\varepsilon,D;M)=\{(K,M):K\in C^{\infty},$ $\dim K=k,$
$K\subset D,$ and $i(K,M)\geq\varepsilon\}$
\end{definition}

\begin{definition}
For $K\in\mathcal{A}^{\infty}(k,\varepsilon,D;M),$

i. $Sect(K)$ denotes the sectional curvatures of $K$ and

ii. $II_{w}^{K}$ denotes the second fundamental form of $K$ with respect to a
normal vector $w$ in $M.$

Define $\left\|  II^{K}\right\|  (p)=\max\{\left\|  II_{w}^{K}(v)\right\|
:\forall w\in UNK_{p}$ and $v\in UTK_{p}\}$ and $\left\|  II^{K}\right\|
=\sup_{p\in K}\left\|  II^{K}\right\|  (p).$
\end{definition}

\begin{proposition}
Let $\mathcal{A}^{\infty}(k,\varepsilon,D;M)$ be given and $g_{0}$ be the
Riemannian metric of $M$. There exists positive constants $C_{0},C_{1}%
,d_{0},v_{0},i_{0}$ depending on $n,k,\varepsilon,D$ and $M$ such that
$\forall K\in\mathcal{A}^{\infty}(k,\varepsilon,D;M),$ $e:K\hookrightarrow M,$
the Riemannian manifold $(K,e^{\ast}g_{0})$ satisfies the following intrinsically:

i. $\left\|  II^{K}\right\|  \leq C_{0}$ and $\left|  Sect(K)\right|  \leq
C_{1}, $

ii. $v(K)\geq v_{0},$

iii. $d(K)\leq d_{0}$, and

iv. consequently, $i(K)\geq i_{0}$ by Cheeger [Ch], [CE].
\end{proposition}

\begin{proof}
Let $D=\overline{B}(D,\varepsilon)$ and $\varepsilon_{0}=\min(\varepsilon
,\frac{1}{2}i(D^{\prime})).$

i. $\forall q\in D^{\prime},$ $\partial B(q,\varepsilon_{0})$ are smooth
submanifolds of $M$, which are diffeomorphic to $S^{n-1}.$ By compactness,
there exists $C_{0}>0$ such that $\left\|  II^{\partial B(q,\varepsilon_{0}%
)}\right\|  \leq C_{0}$, for all $q\in D^{\prime}.$

Let $K\in\mathcal{A}^{\infty}(k,\varepsilon,D;M)$ be arbitrarily chosen.
Choose any $p\in K,$ $v\in UK_{p},$ and $w\in UNK_{p}.$ Let $q=\exp_{p}%
^{M}\varepsilon_{0}w.$ Since $\varepsilon_{0}\leq\varepsilon<R_{O}(K,M),$
$B(q,\varepsilon_{0})\cap K=\emptyset.$ Let $S$ denote $\partial
B(q,\varepsilon_{0})$ below in this proof. $S$ is smooth and $v\in UTS_{p}.$
Define $\alpha_{1}(s)=\exp_{p}^{K}sv$ and $\alpha_{2}(s)=\exp_{p}^{S}sv,$
$f_{j}(s)=d_{M}(\alpha_{j}(s),q),$ for $j=1,2$, and $W=-grad$ $d_{M}(.,q).$
$f_{1}$ has a local minimum at $s=0,$ hence $f_{1}^{\prime\prime}(0)\geq0=$
$f_{2}^{\prime\prime}(0).$%
\begin{align*}
f_{j}^{\prime\prime}(0)  &  =-\left(  \left\langle \nabla_{\alpha_{j}^{\prime
}}^{M}\alpha_{j}^{\prime},w\right\rangle +\left\langle \nabla_{v}%
^{M}W,v\right\rangle \right) \\
\left\langle \nabla_{\alpha_{1}^{\prime}}^{M}\alpha_{1}^{\prime}%
,w\right\rangle  &  \leq\left\langle \nabla_{\alpha_{2}^{\prime}}^{M}%
\alpha_{2}^{\prime},w\right\rangle \\
II_{w}^{K}(v)  &  \leq II_{w}^{S}(v)\\
-II_{w}^{K}(v)  &  =II_{-w}^{K}(v)\leq II_{-w}^{S^{\prime}}(v)
\end{align*}
where $S^{\prime}=\partial B\left(  \exp_{p}^{M}\left(  -\varepsilon
_{0}w\right)  ,\varepsilon_{0}\right)  .$ Hence, $\left\|  II^{K}\right\|
\leq C_{0}.$ By Gauss's Theorem [DoC, p130] relating the second fundamental
form and the sectional curvatures of $K$ and $M,$ and the polarization
identities, there exists $C_{1}(C_{0},\left|  Sect(M)\right|  )$ such that
$\left|  Sect(K)\right|  \leq C_{1}.$

ii. There exists $v_{1}>0$ such that $\forall p\in D^{\prime},vol_{n}%
(B(q,\varepsilon;M))\geq v_{1}$ by compactness of $D^{\prime}.$ Furthermore,
$v_{1}$ can be chosen only depending on the dimension $n,\varepsilon$ and
$i(D^{\prime})$ but not on $M,$ by using the estimates of the lower bounds for
the volumes of the balls of radius less than $i(D^{\prime})/2$ by Croke [Cr,
Prop.14]. Let $K\in\mathcal{A}^{\infty}(k,\varepsilon,D;M)$ and $p\in K$ be
arbitrarily chosen. By \emph{Theorem 2.1} and \emph{Remark 2}, page 453 of
Heintze \& Karcher [HK]:
\[
0<v_{1}\leq vol_{n}(B(p,\varepsilon;M))\leq vol_{n}(B(K,\varepsilon;M))\leq
v(K)\cdot C(k,n,C_{0},C_{1},\varepsilon)
\]
where $v(K)$ is the k-dimensional volume of the $K$ with the induced
submanifold metric.

iii. Choose $d_{2}=\min(\varepsilon_{0},d_{1})$ with $d_{1}$ of Lemma 2 below.
There exists $v_{2}>0$ such that $\forall p\in D^{\prime},vol_{n}(B(q,\frac
{1}{4}d_{2};M))\geq v_{2}$ as in part (ii). Let $K\in\mathcal{A}^{\infty
}(k,\varepsilon,D;M)$ be arbitrarily chosen. Since all geodesics $\gamma$ of
$K$ satisfy $\left\|  \nabla_{\gamma^{\prime}}^{M}\gamma^{\prime}\right\|
\leq C_{0}$ by part (i), one can conclude that
\[
\forall p\in K,\text{ }B(p,\frac{1}{2}d_{2};M)\cap K\subset B(p,d_{2};K)
\]
by using the Lemma 2(i) and $d_{2}\leq i(K,M).$ Let $p$ and $q$ be a pair of
intrinsically furthest apart points in $K,$ that $d_{K}(p,q)=d(K,e^{\ast}%
g_{0}).$ Choose a normal minimal geodesic $\gamma$ of $K$ form $p$ to $q,$
$\gamma(0)=p,\gamma(d(K))=q,$ and $\left\|  \gamma^{\prime}\right\|  =1.$
Suppose that $B(\gamma(ad_{2}),\frac{1}{4}d_{2};M)\cap B(\gamma(bd_{2}%
),\frac{1}{4}d_{2};M)\neq\emptyset,$ for some integers $a,b\in\mathbf{N}%
\cap\lbrack0,\frac{d(K)}{d_{2}}].$ Then $d_{M}(\gamma(ad_{2}),\gamma
(bd_{2}))<\frac{1}{2}d_{2}$ which implies $d_{K}(\gamma(ad_{2}),\gamma
(bd_{2}))<d_{2}.$ Thus, $a=b.$ Hence, the balls $B(\gamma(ad_{2}),\frac{1}%
{4}d_{2};M),$ for $a\in\mathbf{N}\cap\lbrack0,\frac{d(K)}{d_{2}}],$ are
disjoint in $D^{\prime}.$%
\[
v_{2}\cdot\frac{d(K)}{d_{2}}\leq vol_{n}(D^{\prime}).
\]

iv. This follows Cheeger [Ch] and parts (i-iii).
\end{proof}

\begin{lemma}
Let $D^{\prime}$ be a compact subset of $M$. Given $C_{0},$ there exist
$0<d_{1}\leq\frac{1}{2}i(D^{\prime})$ such that any $C^{2}$ curve
$\alpha:[0,d_{1}]\rightarrow D^{\prime}$ with $\left\|  \alpha^{\prime
}(s)\right\|  =1$ and $\left\|  \nabla_{\alpha^{\prime}}\alpha^{\prime
}\right\|  \leq C_{0},$ must satisfy

i. $d_{M}(\alpha(0),\alpha(s))\geq\frac{3s}{4},\forall s\in\lbrack0,d_{1}]$ and

ii. $d_{M}(\gamma(s),\alpha(s))\leq\frac{s}{4},\forall s\in\lbrack0,d_{1}]$
where $\gamma(s)=\exp_{\alpha(0)}s\alpha^{\prime}(0)$ $.$
\end{lemma}

\begin{proof}
i. Suppose that such $d_{1}>0$ does not exist. Then, $\forall m\in
\mathbf{N}^{+},\exists\alpha_{m}:[0,1]\rightarrow D^{\prime}$ with
$d(\alpha_{m}(0),\alpha_{m}(s_{m}))<\frac{3s_{m}}{4}$ for some $s_{m}%
\in(0,\frac{1}{m}],$ $\left\|  \alpha_{m}^{\prime}(s)\right\|  =1$ and
$\left\|  \nabla_{\alpha_{m}^{\prime}}\alpha_{m}^{\prime}\right\|  \leq
C_{0}.$ Then by compactness of $D^{\prime}$, there exists a subsequence which
we denote with the same subindices, $\alpha_{m}(0)\rightarrow p_{0}$ and
$\alpha_{m}^{\prime}(0)\rightarrow v_{0}\in UTM_{p_{0}}.$ Since $d\exp_{p_{0}%
}(0)=Id$, for a given $\delta>0,$ there are sufficiently small $\eta>0$,
$\sigma>0$ and sufficiently large $m$ such that $\tilde{\alpha}_{m}(s)=\left(
\exp_{p_{0}}|B(0,\eta,TM_{p})\right)  ^{-1}\alpha_{m}(s)$ are defined for
$0\leq s\leq\sigma$, $\left|  \left\|  \tilde{\alpha}_{m}^{\prime}(s)\right\|
-1\right|  <\delta,$ $\left\|  \tilde{\alpha}_{m}(0)-\tilde{\alpha}_{m}%
(s_{m})\right\|  <\frac{3s_{m}}{4}(1+\delta),$ and $\left\|  \tilde{\alpha
}_{m}^{\prime\prime}(s)\right\|  \leq C_{2}$. $\eta,\sigma$ and $C_{2}$ depend
on $\delta,C_{0}$, the metric $g_{0}$ locally and the derivatives of
$\exp_{p_{0}}$ near $0,$ but not on $m.$ This contradicts Schur's Theorem in
$\mathbf{R}^{n}$, [Cn] or the fact that all $C^{2}$ curves $\gamma$ in
$\mathbf{R}^{n}$, with $\left\|  \gamma^{\prime}(s)\right\|  =1$ and $\left\|
\gamma^{\prime\prime}(s)\right\|  \leq C_{2}$ satisfy $\left\|  \gamma
(s)-\gamma(0)\right\|  \geq\frac{\sin sC_{2}}{C_{2}} $ for $s\in(0,\frac{\pi
}{2C_{2}}]$, (see [D6, proof of Proposition 2a] for a proof). Consequently,
$\exists d_{1}>0$ as indicated.

ii. This is an immediate consequence of (i) since $\gamma(0)=\alpha(0)$ and
$d(\gamma(s),\gamma(0))=s$ for $0\leq s\leq d_{1}\leq i(D^{\prime}).$
\end{proof}

\begin{proposition}
Let $\mathcal{A}^{\infty}(k,\varepsilon,D;M)$ be given and $g_{0}$ be the
Riemannian metric of $M$. Consider
\[
\mathcal{C}^{\infty}(k,\varepsilon,D;M)=\{(K,e^{\ast}g_{0}):\forall
K\in\mathcal{A}^{\infty}(k,\varepsilon,D;M),e:K\hookrightarrow D\subset M\}
\]
as a collection of Riemannian manifolds, not as submanifolds of $M.$ By
Gromov's (pre)Compactness Theorem, $\mathcal{C}^{\infty}(k,\varepsilon,D;M) $
has finitely many diffeomorphism types, and any sequence $(K_{m},g_{m})$ in
$\mathcal{C}^{\infty}(k,\varepsilon,D;M)$ has a Cauchy subsequence $(K_{m_{j}%
},g_{m_{j}})$ in $\mathcal{C}^{\infty}(k,\varepsilon,D;M)$ with respect to
Lipschitz distance. As it was stated in [Pe], all $K_{m_{j}}$ are
diffeomorphic to a fixed $C^{\infty}$ manifold $K,$ and $g_{m_{j}}\rightarrow
g_{\infty}$ in $C^{1}$ sense on $K$ with respect to some harmonic
coordinates$,$ in which $g_{\infty}$ is a $C^{1,\alpha}$ Riemannian metric of
$K.$ $(K,g_{\infty})$ is an Alexandrov space of bounded curvature by [Ni].
\end{proposition}

\begin{proof}
In Proposition 8, we proved all necessary conditions for hypothesis of
Gromov's Compactness Theorem, see [Gr], [Ni], [Pe], [GW] and [D3].
\end{proof}

\begin{remark}
A priori, $(K,g_{\infty})$ is not a Riemannian submanifold of $(M,g_{0}).$ We
will prove in Theorem 2 that there exists an isometric embedding
$(K,g_{\infty})\hookrightarrow(M,g_{0}).$ If one starts with an arbitrary
$C^{1,1}$ $K\in\mathcal{A}(k,\varepsilon,D;M),$ $e:K\hookrightarrow D\subset
M,$ then $e^{\ast}g_{0}$ is a priori $C^{0,1}$, which is of too low regularity
to have any sense of curvature. It is not a priori necessary that smooth
approximations of $e^{\ast}g_{0}$ are of uniformly bounded curvature or smooth
approximations of the embedding $e:K\hookrightarrow D\subset M$ have thickness
close to $\varepsilon.$
\end{remark}

\begin{notation}
For $C^{1}$ $f:\mathbf{R}^{k}\rightarrow\mathbf{R}^{m},v\in U\mathbf{R}^{k}$
the directional derivative of $f$ in the direction $v$ is $f_{v}(p)=\left.
\frac{d}{dt}f(p+tv)\right|  _{t=0}.$ The Jacobian $f^{\prime}(p)$ is an
$m\times k$ matrix and $\left\|  f^{\prime}(p)\right\|  $ is its norm in
$\mathbf{R}^{km}$, and $\left\|  f^{\prime}\right\|  =\underset{p}{\sup
}\left\|  f^{\prime}(p)\right\|  .$ For $p\in\mathbf{R}^{k}$, $v\in
U\mathbf{R}_{p}^{k}$, $f_{vv}(p)=\left.  \frac{d^{2}}{dt^{2}}f(p+tv)\right|
_{t=0}$ which is defined almost everywhere in $p$, when $f\in C^{1,1}.$
\end{notation}

\begin{lemma}
Let $f:\mathbf{R}^{k}\rightarrow\mathbf{R}^{m}$ be $C^{1,1}$, $p\in
\mathbf{R}^{k}$, $v\in U\mathbf{R}_{p}^{k}$ and let $G$ be the graph of $f$ in
$\mathbf{R}^{k+m}$. Define
\begin{align*}
I(f,p,v,w)  &  =\left(  1+\left\|  f_{v}(p)\right\|  ^{2}\right)  \left(
1+\left\|  \nabla(f\cdot w)(p)\right\|  ^{2}\right)  ^{\frac{1}{2}}\geq1,\\
\forall x  &  \in\mathbf{R}^{k},\forall v\in U\mathbf{R}^{k},w\in
U\mathbf{R}^{n-k}.
\end{align*}
If $F_{g}((p,f(p)),G)\geq R$ and $f_{vv}(p)$ exists, then
\[
\forall w\in U\mathbf{R}^{m},\left\|  f_{vv}(p)\cdot w\right\|  \leq\frac
{1}{R}I(f,p,v,w).
\]
Conversely, if $f_{vv}(p)$ exists and $\exists w\in U\mathbf{R}^{m}$ such that
$\left\|  f_{vv}(p)\cdot w\right\|  >\frac{1}{R}I(f,p,v,w),$ then
$F_{g}((p,f(p)),G)<R$ and particularly, $B((p,f(p))+Rn,R)\cap G\neq\emptyset$
where $n=(-\nabla(f\cdot w)(p),w)\left(  1+\left\|  \nabla(f\cdot
w)(p)\right\|  ^{2}\right)  ^{-\frac{1}{2}}.$
\end{lemma}

\begin{proof}
Let $\phi(u)=(u,f(u))=(x,y)\in\mathbf{R}^{k}\times\mathbf{R}^{m}$ be a
parametrization of $G.$ For $(w_{1},w)\in\mathbf{R}^{k}\times\mathbf{R}^{m} $
to be normal to $G$ at $p,$ $(u,f_{u}(p))\cdot(w_{1},w)=0$ should be true for
all $u\in U\mathbf{R}_{p}^{k}$, that is $(w_{1},w)\in nullspace([I_{k}$
$|f^{\prime}(p)^{T}]).$%
\[
(w_{1},w)=\sum\nolimits_{j=1}^{m}(-\nabla f_{j}(p),e_{j})w_{j}=(-\nabla(f\cdot
w)(p),w)
\]
Let $n=(w_{1},w)\left\|  (w_{1},w)\right\|  ^{-1}=(w_{1},w)\left(  \left\|
w\right\|  ^{2}+\left\|  \nabla(f\cdot w)(p)\right\|  ^{2}\right)  ^{-\frac
{1}{2}}$

and define $\sigma(t)=\frac{1}{2}\left\|  \phi(p)+Rn-\phi(p+tv)\right\|
^{2}.$%
\begin{align*}
\sigma^{\prime}(0)  &  =-(v,f_{v}(p))\cdot Rn=0\\
\sigma^{\prime\prime}(0)  &  =-(0,f_{vv}(p))\cdot Rn+\left\|  (v,f_{v}%
(p))\right\|  ^{2}\\
\text{\ }f_{vv}(p)\cdot w  &  =\frac{1}{R}\left(  1+\left\|  f_{v}(p)\right\|
^{2}-\sigma^{\prime\prime}(0)\right)  \left(  1+\left\|  \nabla(f\cdot
w)(p)\right\|  ^{2}\right)  ^{\frac{1}{2}},\forall w\in U\mathbf{R}^{m}%
\end{align*}
If $F_{g}(p,G)\geq R$, then $B(\phi(p)+Rn,R)\cap G=\emptyset$ and $\sigma(t)$
has a local minimum at $t=0,\,$that is $\sigma^{\prime\prime}(0)\geq0,$ since
it exists.
\[
\text{\ }f_{vv}(p)\cdot w\leq\frac{1}{R}\left(  1+\left\|  f_{v}(p)\right\|
^{2}\right)  \left(  1+\left\|  \nabla(f\cdot w)(p)\right\|  ^{2}\right)
^{\frac{1}{2}},\forall w\in U\mathbf{R}^{m}%
\]
Using $-w$ gives the inequality with the absolute value$.$ For the converse,
choose $w$ or $-w$ for positive $f_{vv}(p)\cdot w.$ $\sigma^{\prime\prime
}(0)<0$ implies $\sigma(t)<\sigma(0)$ for small $t\neq0,$ even for a $C^{1,1}$
function. Hence, $B(\phi(p)+Rn,R)\cap G\neq\emptyset$ and $F_{g}(p,G)<R.$
\end{proof}

\begin{lemma}
Let $f:\mathbf{R}^{k}\rightarrow\mathbf{R}^{m}$ be $C^{1,1}$ and satisfy

\qquad a. $\left\|  f^{\prime}\right\|  \leq A,$

\qquad b. $\left\|  f^{\prime}(x)-f^{\prime}(y)\right\|  \leq B\left\|
x-y\right\|  ,\forall x,y\in\mathbf{R}^{k},$ and

\qquad c.$\left\|  f_{vv}(x)\cdot w\right\|  \leq C$ for a.e. $x\in
\mathbf{R}^{k},$ for fixed $v\in U\mathbf{R}^{k}$ and $w\in U\mathbf{R}^{m}.$

Then $\forall\delta,\rho>0,\exists$ a $C^{1,1}$ function $h:\mathbf{R}%
^{k}\rightarrow\mathbf{R}^{m}$ such that

\qquad i. $h$ is $C^{\infty}$ on $B(0,\rho)$ and $h=f$ outside $B(0,2\rho),$

\qquad ii. $\left\|  h-f\right\|  \leq A\delta,$

\qquad iii. $\left\|  h^{\prime}-f^{\prime}\right\|  \leq(aA+B)\delta,$

\qquad iv. $\left\|  h^{\prime}(x)-h^{\prime}(y)\right\|  \leq\left\|
x-y\right\|  (B+\delta(2aB+bA)),$ and

\qquad v. $\left\|  h_{vv}(x)\cdot w\right\|  \leq C+\delta(2aB+bA)).$

where $a$ and $b$ are constants depending on $\frac{1}{\rho}$ but not on $f. $
\end{lemma}

\begin{proof}
Choose $\eta:[0,\infty)\rightarrow\lbrack0,1]$ smooth with $\sup
p(\eta)\subset\lbrack0,1],\eta^{-1}(1)=[0,\frac{1}{2}],$

$-2.25\leq\eta^{\prime}\leq0,$ and $\left\|  \eta^{\prime\prime}\right\|  \leq10.$

Define $g:\mathbf{R}^{k}\rightarrow\mathbf{R}^{m}$ and $h:\mathbf{R}%
^{k}\rightarrow\mathbf{R}^{m}$ by
\begin{align*}
g(x)  &  =c_{\delta}\int_{\left\|  u\right\|  \leq\delta}f(x+u)\eta\left(
\frac{\left\|  u\right\|  }{\delta}\right)  du\text{ where }c_{\delta}%
^{-1}=\int_{\left\|  u\right\|  \leq\delta}\eta\left(  \frac{\left\|
u\right\|  }{\delta}\right)  du\\
h(x)  &  =(1-\eta(2\rho\left\|  x\right\|  ))f+\eta(2\rho\left\|  x\right\|
)g.
\end{align*}

Set $a=\sup\left\|  \left(  \eta(2\rho\left\|  x\right\|  )\right)  ^{\prime
}\right\|  $ and $b=\sup\left\|  \left(  \eta(2\rho\left\|  x\right\|
)\right)  ^{\prime\prime}\right\|  .$

The proofs of (i-iv) are elementary and will be left to the reader. We will
only give a proof of (v). Let $\lambda(t)=f(x+tv)\cdot w.$ Then $\lambda
^{\prime}(t)=f_{v}(x+tv)\cdot w$ which is lipschitz and hence absolutely
continuous, and $\lambda^{\prime\prime}(t)=f_{vv}(x+tv)\cdot w$ which is
defined almost everywhere.
\begin{align*}
\left\|  \left(  f_{v}(x+tv)-f_{v}(x)\right)  \cdot w\right\|   &  =\left\|
\lambda^{\prime}(t)-\lambda^{\prime}(0)\right\|  =\left\|  \int\limits_{0}%
^{t}\lambda^{\prime\prime}(u)du\right\|  \leq Ct,\forall x\\
\left\|  \left(  g_{v}(x+tv)-g_{v}(x)\right)  \cdot w\right\|   &  \leq
c_{\delta}\int_{\left\|  u\right\|  \leq\delta}\left\|  \left(  f_{v}%
(x+tv+u)-f_{v}(x+u)\right)  \cdot w\right\|  \eta\left(  \frac{\left\|
u\right\|  }{\delta}\right)  du\\
&  \leq Ct,\forall x\\
\left\|  g_{vv}(x)\cdot w\right\|   &  \leq C,\forall x\text{ a.e.}\\
\left\|  (g-f)(x)\cdot w\right\|   &  \leq c_{\delta}\int_{\left\|  u\right\|
\leq\delta}\left\|  \left(  f(x+u)-f(x)\right)  \cdot w\right\|  \eta\left(
\frac{\left\|  u\right\|  }{\delta}\right)  du\\
&  \leq A\delta,\forall x,\\
\left\|  (g_{v}-f_{v})(x)\cdot w\right\|   &  \leq c_{\delta}\int_{\left\|
u\right\|  \leq\delta}\left\|  \left(  f_{v}(x+u)-f_{v}(x)\right)  \cdot
w\right\|  \eta\left(  \frac{\left\|  u\right\|  }{\delta}\right)  du\\
&  \leq B\delta,\forall x\\
\forall x\text{ a.e.},\left\|  h_{vv}(x)\cdot w\right\|   &  \leq\left\|
(1-\eta(2\rho\left\|  x\right\|  ))f_{vv}\cdot w+\eta(2\rho\left\|  x\right\|
)g_{vv}\cdot w\right\|  +\\
&  2\left\|  \eta(2\rho\left\|  x\right\|  )_{v}(g-f)_{v}\cdot w\right\|
+\left\|  \eta(2\rho\left\|  x\right\|  )_{vv}(g-f)\cdot w\right\| \\
& \\
&  \leq C+2aB\delta+bA\delta
\end{align*}
\end{proof}

\begin{proposition}
i. Let $K$ be a compact $C^{1,1}$ submanifold of $\mathbf{R}^{n}$ such that
$F_{g}(K)\geq R_{1}.$ Then $\forall R_{2}<R_{1},\forall\sigma>0,$ there exists
a smooth approximation $K^{\delta}$ of $K$ in $M$ such that $d_{C^{1}%
}(K,K^{\delta})<\sigma$ and $F_{g}(K^{\delta})\geq R_{2}.$

ii. Hence, $\forall\sigma>0,$ $\mathcal{A}(k,\varepsilon,D;\mathbf{R}%
^{n})\subset\overline{\mathcal{A}^{\infty}(k,\varepsilon,D_{\sigma}%
;\mathbf{R}^{n})}$ with respect to $C^{1}$ topology, where $D_{\sigma}$ is the
closure of the $\sigma-$neighborhood of $D.$
\end{proposition}

\begin{proof}
i. Let $p\in K$ be any point. Rotate and translate $K$ in $\mathbf{R}%
^{k}\times\mathbf{R}^{n-k}$ so that $p=0$ and $TK_{p}=\mathbf{R}^{k}%
\times\mathbf{\{0\}.}$ $\exists r_{p}>0,\exists f_{p}:B(0,3r,\mathbf{R}%
^{k})\rightarrow\mathbf{R}^{n-k}$ such that $U=\{(x,f_{p}(x)):\left\|
x\right\|  <3r_{p}\}$ is an open neighborhood of $p$ in $K,$ and $\left\|
f_{p}^{\prime}(x)\right\|  \leq1$ for $\left\|  x\right\|  <3r.$ Set
$V(p)=\{(x,f_{p}(x)):\left\|  x\right\|  <r_{p}\}.$ By the compactness of $K,$
there are finitely many $V(p_{j})$ covering $K.$ Let $r_{j},f_{j},U_{j},V_{j}
$ be defined as above associated to $p_{j}.$ Choose $k_{j}$ and $A_{j}$ so
that
\[
\frac{1}{R_{1}}=k_{1}<k_{2}<...<k_{j_{0}+1}=\frac{1}{R_{2}}\text{ and }%
1=A_{1}<A_{2}<...<A_{j_{0}+1}=2.
\]
Rotate $K$ so that $p_{1}=0$ and $TK_{p_{1}}=\mathbf{R}^{k}\times
\mathbf{\{0\}.}$ $F_{g}(U_{1})\geq R_{1}$ by the hypothesis. By Lemma 3:
\[
\left\|  f_{1,vv}(x)\cdot w\right\|  \leq k_{1}I(f_{1},x,v,w),\forall x\in
B(0,3r_{1}),\forall v\in U\mathbf{R}^{k},w\in U\mathbf{R}^{n-k},a.e.
\]
By Lemma 4, for a given $\delta>0,$ there exists a $C^{1,1}$ approximation
$h_{1}^{\delta}$ of $f_{1}$ which is $C^{\infty}$ for $\left\|  x\right\|
<r_{1}$ and coincides with $f_{1}$ for $\left\|  x\right\|  \geq2r_{1}.$ Let
$K_{1}^{\delta}$ be the submanifold of $\mathbf{R}^{n}$ obtained from $K$ by
replacing $U_{1}$ with the graph $U_{1}^{\delta}$ of $h_{1}^{\delta}.$
\begin{align*}
\left\|  f_{1}-h_{1}^{\delta}\right\|   &  \leq A_{1}\delta\\
\left\|  \left(  f_{1}\right)  ^{\prime}-\left(  h_{1}^{\delta}\right)
^{\prime}\right\|   &  \leq\delta\cdot c_{3}(A_{1},R_{1},\left\|
f_{1}^{\prime}\right\|  )\\
\text{ and }\forall x  &  \in B(0,3r_{1}),\forall v\in U\mathbf{R}^{k},w\in
U\mathbf{R}^{n-k},a.e.,\\
\left\|  h_{1,vv}^{\delta}(x)\cdot w\right\|   &  \leq k_{1}I(f_{1}%
,x,v,w)+\delta\cdot c_{4}(A_{1},R_{1},\left\|  f_{1}^{\prime}\right\|  )
\end{align*}
Choose $\delta_{1}>0$ small enough so that $\forall\delta$ with $0<\delta
\leq\delta_{1},$

a. $(1+\delta)A_{1}\leq A_{2},$ and

b. $\forall x\in B(0,3r_{1}),\forall v\in U\mathbf{R}^{k},w\in U\mathbf{R}%
^{n-k},a.e.,$%

\[
\frac{\left\|  h_{1,vv}^{\delta}(x)\cdot w\right\|  }{I(h_{1}^{\delta}%
,x,v,w)}\leq\frac{k_{1}I(f_{1},x,v,w)+\delta\cdot c_{4}(A_{1},R_{1},\left\|
f_{1}^{\prime}\right\|  )}{I(h_{1}^{\delta},x,v,w)}\leq k_{2},\text{ and }%
\]

c. the adjustment of $f_{1}$ by $h_{1}^{\delta}$ does not change $f_{j}$ being
graphs and keeps $\left\|  f_{j}^{^{\prime}}\right\|  \leq A_{2},$ for $j\geq2.$

Consequently, $F_{g}(K_{1}^{\delta})\geq\frac{1}{k_{2}}$ and $K_{1}^{\delta}$
is smooth on $U_{1}^{\delta}$. One proceeds inductively to obtain a
$C^{\infty}$ approximation $K^{\delta}$ of $K,$ for all $0<\delta\leq
\delta_{j_{0}}$ for some $\delta_{j_{0}}>0.$ Then,

a. $F_{g}(K^{\delta})\geq R_{2},\forall0<\delta\leq\delta_{j_{0}},$ and

b. $\lim_{\delta\rightarrow0}d_{C^{1}}(K,K^{\delta})=0,$ that is
$\forall\sigma>0,\exists K^{\delta}$ such that $d_{C^{1}}(K,K^{\delta})<\sigma.$

\smallskip

ii. Let $\sigma>0$ and $K\in\mathcal{A}(k,\varepsilon,D;\mathbf{R}^{n})$ be
given, that is $F_{g}(K)\geq i(K,\mathbf{R}^{n})=R_{O}(K,\mathbf{R}%
^{n})>\varepsilon.$ Then, $\forall m\in\mathbf{N}^{+}$ with $m>\frac{1}%
{\sigma},\mathbf{\exists}K_{m},$ a smooth approximation of $K$ in
$\mathbf{R}^{n}$ such that $F_{g}(K)-\frac{1}{m}\leq F_{g}(K_{m})$ and
$d_{C^{1}}(K,K_{m})\leq\frac{1}{m}.$

Claim 1. $\lim\inf_{m}MDC(K_{m})>0$. Suppose that $\lim\inf_{m}MDC(K_{m})=0 $
and follow the proof of Proposition 3, to obtain $p_{0}=q_{0}.$ If
$\eta_{p_{m}q_{m}}$ denotes a minimal geodesic of $K_{m}$ between $p_{m}$ and
$q_{m},$ then $\eta_{p_{m}q_{m}}$ is normal to the segment $\gamma_{p_{m}%
q_{m}}$ at $p_{m}$ and $q_{m}.$ Since $\left\|  p_{m}-q_{m}\right\|
\rightarrow0,$ the maximum of the ambient curvature of $\eta_{p_{m}q_{m}}$ in
$\mathbf{R}^{n}$ becomes arbitrarily large as $m\rightarrow\infty.$ But, the
sectional curvatures of $K_{m}$ are bounded by $\frac{2}{\varepsilon}$ for
large $m$ by Proposition 8(i). Thus, Claim 1 holds.

By Propositions 2 and 3:
\[
\underset{m\rightarrow\infty}{\lim\sup}R_{O}(K_{m},\mathbf{R}^{n})\leq
R_{O}(K,\mathbf{R}^{n})\leq\frac{1}{2}MDC(K,\mathbf{R}^{n})\leq\underset
{m\rightarrow\infty}{\lim\inf}\frac{1}{2}MDC(K_{m},\mathbf{R}^{n}).
\]

Claim 2. $\underset{m\rightarrow\infty}{\lim\sup}R_{O}(K_{m},\mathbf{R}%
^{n})=R_{O}(K,\mathbf{R}^{n}).$

Suppose\ that $\underset{}{\underset{m\rightarrow\infty}{\lim\sup}}R_{O}%
(K_{m},\mathbf{R}^{n})<R_{O}(K,\mathbf{R}^{n}).$ Then for sufficiently large
$m,$
\[
R_{O}(K_{m},\mathbf{R}^{n})<\frac{1}{2}MDC(K_{m},\mathbf{R}^{n}),\text{i.e.
}R_{O}(K_{m},\mathbf{R}^{n})=F_{g}(K_{m},\mathbf{R}^{n}).
\]
However, this brings all to a contradiction:
\[
R_{O}(K,\mathbf{R}^{n})\leq F_{g}(K,\mathbf{R}^{n})\underset{}{\leq
\underset{m\rightarrow\infty}{\lim\sup}F_{g}(K_{m},\mathbf{R}^{n}%
)=\underset{m\rightarrow\infty}{\lim\sup}}R_{O}(K_{m},\mathbf{R}^{n}%
)<R_{O}(K,\mathbf{R}^{n})
\]
Hence, $\underset{m\rightarrow\infty}{\lim\sup}R_{O}(K_{m},\mathbf{R}%
^{n})=R_{O}(K,\mathbf{R}^{n})>\varepsilon,$ where the smooth submanifolds
$K_{m}\subset D_{\sigma}$ and $K_{m}\rightarrow K$ in $C^{1}$ sense. In other
words, $K\in\overline{\mathcal{A}^{\infty}(k,\varepsilon,D_{\sigma}%
;\mathbf{R}^{n})}.$
\end{proof}

\begin{proposition}
i. For any given $\varepsilon>0,$ a complete Riemannian manifold $M$ and a
compact subset $D\subset M,$ there exists $\varepsilon^{\prime}(\varepsilon
,D_{\varepsilon},M)>0$ with $\varepsilon^{\prime}<\varepsilon$ satisfying that

''$\forall\sigma>0$ and for any given compact $C^{1,1}$ submanifold $K$ of $M$
with $K\subset D$ and $F_{g}(K)>\varepsilon$, there exists a smooth
approximation $K^{\prime}$ of $K$ in $M$ with $d_{C^{1}}(K,K^{\prime})<\sigma$
and $F_{g}(K^{\prime})>\varepsilon^{\prime}$''.

ii. Hence, $\forall\sigma>0,$ $\mathcal{A}(k,\varepsilon,D;M)\subset
\overline{\mathcal{A}^{\infty}(k,\varepsilon^{\prime},D_{\sigma};M)}$ with
respect to $C^{1}$ topology, where $D_{\sigma}$ is the closure of the
$\sigma-$neighborhood of $D.$
\end{proposition}

\begin{proof}
i. Let $D^{\prime}=\overline{B}(D,\varepsilon)$ and $r_{0}=\frac{1}%
{4}i(D^{\prime})>0.$ Choose a finite collection of points $p_{\alpha}$ such
that $\{B(p_{\alpha},r_{0};M):\alpha=1,...,\alpha_{0}\}$ covers $D^{\prime}$
and let $\varphi_{\alpha}:=\left(  \exp_{p_{\alpha}}^{M}|B(0,3r_{0}%
;TM_{p_{\alpha}})\right)  ^{-1}.$

Define $\varepsilon_{0}=\min(\varepsilon,\frac{1}{2}i(D^{\prime}))$ and
$S(p,\alpha,\varepsilon_{0})=\varphi_{\alpha}(B(p_{\alpha},2r_{0}%
;M)\cap\partial B(p,\varepsilon_{0};M))$ which are smooth since $\varepsilon
_{0}\leq\frac{1}{2}i(D^{\prime}).$ In all of the second fundamental form
assertions below, $\partial B(p,\varepsilon_{0};M)$ are codimension 1 smooth
submanifolds of $(M,g_{0}),$ and $S(p,\alpha,\varepsilon_{0})$ are codimension
1 smooth submanifolds of $\mathbf{R}^{n}$ with the flat metric.
\begin{align*}
\exists c_{5},\forall p  &  \in D^{\prime},\left\|  II^{\partial
B(p,\varepsilon_{0};M)}\right\|  \leq c_{5}\\
\exists c_{6},\forall p  &  \in D^{\prime},\forall\alpha,\left\|
II^{S(p,\alpha,\varepsilon_{0})}\right\|  \leq c_{6}\text{ whenever
}S(p,\alpha,\varepsilon_{0})\neq\emptyset
\end{align*}

The first assertion follows the smoothness of the metric $g_{0}$ of $M,$
$\varepsilon_{0}\leq\frac{1}{2}i(D^{\prime}),$ and compactness of $D^{\prime
}.$ The second assertion follows the facts that there are finitely many
$\alpha,$ and $\left(  \varphi_{\alpha}\right)  _{\ast}g_{0}$ are uniformly
$C^{i}-$bounded for $i=0,1,2,$ as well as quasi-isometric to the Euclidean
metric: $\infty>a\geq\frac{\left\|  \left(  \varphi_{\alpha}\right)  _{\ast
}g_{0}(v)\right\|  }{\left\|  v\right\|  }\geq b>0,\ $uniformly. $\forall
\tau>0,$ define
\begin{align*}
\underset{p,v,w}{\min}\left\|  II_{w}^{\partial B(p,\tau;M)}(v)\right\|   &
=\lambda(\tau,D^{\prime})\geq0\\
\underset{p,\alpha,v^{\prime},w^{\prime}}{\min}\left\|  II_{w}^{S(p,\alpha
,\tau)}(v)\right\|   &  =\mu(\tau,D^{\prime})\geq0
\end{align*}
where $p\in D^{\prime},q\in\partial B(p,\tau;M),$ $v\in UT\partial
B(p,\tau;M)_{q}$, $w\in UN\partial B(p,\tau;M)_{q}$, $\alpha=1,...,\alpha
_{0},$ $q^{\prime}\in S(p,\alpha,\tau)\neq\emptyset$, $v^{\prime}\in
UT\partial B(p,\tau;M)_{q^{\prime}}$, and $w^{\prime}\in UN\partial
B(p,\tau;M)_{q^{\prime}}.$ Then $\lim_{\tau\rightarrow0^{+}}\lambda
(\tau,D^{\prime})=\infty,$ since $D^{\prime}$ is compact. By the reasons
stated above, also $\mu(\tau,D^{\prime})\rightarrow\infty$, as $\tau
\rightarrow0^{+}.$

Choose any $c_{7}>c_{6}$, $\varepsilon_{1}>0$ and $\varepsilon^{\prime
}(\varepsilon,D^{\prime},M)>0$ such that $\mu(\varepsilon_{1},D^{\prime
})>c_{7}$ and $\varepsilon_{1}>\varepsilon^{\prime}.$

Let $K$ be any given compact $C^{1,1}$ submanifold of $M$ with $K\subset D$
and $F_{g}(K)>\varepsilon.$ Then, $K$ avoids tangential balls $B(p,\varepsilon
_{0},M)$ in an open neighborhood $W$ of the point of tangency $q\in\partial
B(p,\varepsilon_{0},M)$. Then any nonempty $\varphi_{\alpha}(K)$ avoids the
open set $\varphi_{\alpha}(B(p,\varepsilon_{0},M))$ in an open neighborhood
$W^{\prime}$ of the point of tangency $\varphi_{\alpha}(q)\in S(p,\alpha
,\varepsilon_{0})\subset\partial\varphi_{\alpha}(B(p,\varepsilon_{0},M)).$ Let
$n$ be the unit normal to $S(p,\alpha,\varepsilon_{0})$ at $\varphi_{\alpha
}(q)$ towards $\varphi_{\alpha}(B(p,\varepsilon_{0},M)).$ Then for any
$r<\frac{1}{c_{6}},$ $W^{\prime}\cap B(\varphi_{\alpha}(q),r,\mathbf{R}%
^{n})\subset W^{\prime}\cap\varphi_{\alpha}(B(p,\varepsilon_{0},M))$ and
$W^{\prime}\cap\varphi_{\alpha}(K)\cap B(\varphi_{\alpha}(q),r,\mathbf{R}%
^{n})=\emptyset$. Hence,
\[
F_{g}(\varphi_{\alpha}(K\cap B(p_{\alpha},2r_{0};M)),\mathbf{R}^{n})\geq
\frac{1}{c_{6}},\forall\alpha.
\]
By applying the method of Proposition 10 to $\varphi_{\alpha}(K\cap
\overline{B(p_{\alpha},2r_{0};M)}),$ for any $c$ with $c_{6}<c<c_{7},$ there
exists a $C^{1,1}$ approximation $K_{\alpha}$ of $K$ such that $\varphi
_{\alpha}(K_{\alpha}\cap B(p_{\alpha},r_{0};M))$ is a smooth submanifold of
$\mathbf{R}^{n}$, $F_{g}(\varphi_{\alpha}(K_{\alpha}\cap B(p_{\alpha}%
,2r_{0};M)),\mathbf{R}^{n})\geq\frac{1}{c},$ $K$ coincides with $K_{\alpha}$
outside $B(p_{\alpha},3r_{0};M)$ and $d_{C^{1}}(K,K_{\alpha})<\frac{\sigma
}{\alpha_{0}}.$ One proceeds inductively on finitely many $\alpha$ to obtain a
$C^{\infty}$ submanifold $K^{\prime}$ of $M$ satisfying
\begin{align*}
d_{C^{1}}(K,K^{\prime})  &  <\sigma\text{,}\\
F_{g}(\varphi_{\alpha}(K^{\prime}\cap B(p_{\alpha},2r_{0};M)),\mathbf{R}^{n})
&  >\frac{1}{c_{7}},\forall\alpha,\text{ and}\\
\left\|  II^{\varphi_{\alpha}(K^{\prime}\cap B(p_{\alpha},2r_{0};M)}\right\|
&  <c_{7}%
\end{align*}
$\varphi_{\alpha}(K^{\prime}\cap B(p_{\alpha},2r_{0};M))$ avoids tangential
submanifolds $S(p,\alpha,\varepsilon_{1})$ in a deleted open neighborhood of
$\varphi_{\alpha}(q)$ since $\left\|  II_{w}^{S(p,\alpha,\varepsilon_{1}%
)}(v)\right\|  >c_{7}$ for all possible choices of $p,v$ and $w.$ Hence,
$F_{g}(K^{\prime},M)\geq\varepsilon_{1}>\varepsilon^{\prime}.$

ii. Let $\sigma>0$ and $K\in\mathcal{A}(k,\varepsilon,D;M)$ be given:
$F_{g}(K)\geq i(K,M)=R_{O}(K,M)>\varepsilon.$ Then, $\forall m\in
\mathbf{N}^{+}$\textbf{\ }with $m>\frac{1}{\sigma},\mathbf{\exists}K_{m},$ a
smooth approximation of $K$ in $M$ such that $\varepsilon_{1}\leq F_{g}%
(K_{m},M)$ and $d_{C^{1}}(K,K_{m})\leq\frac{1}{m}<\sigma.$

As in Proposition 10, $\lim\inf_{m}MDC(K_{m})>0$.
\[
\underset{m\rightarrow\infty}{\lim\sup}R_{O}(K_{m},M)\leq R_{O}(K,M)\leq
\frac{1}{2}MDC(K,M)\leq\underset{m\rightarrow\infty}{\lim\inf}\frac{1}%
{2}MDC(K_{m},M)
\]
If $\underset{m\rightarrow\infty}{\lim\sup}R_{O}(K_{m},M)=R_{O}%
(K,M)>\varepsilon,$ where the smooth submanifolds $K_{m}\subset D_{\sigma}$
and $K_{m}\rightarrow K$ in $C^{1}$ sense, then,
\[
K\in\overline{\mathcal{A}^{\infty}(k,\varepsilon,D_{\sigma};M)}\subset
\overline{\mathcal{A}^{\infty}(k,\varepsilon^{\prime},D_{\sigma};M)}.
\]
If $\underset{}{\underset{m\rightarrow\infty}{\lim\sup}}R_{O}(K_{m}%
,M)<R_{O}(K,M),$ then for sufficiently large $m$ of the last subsequence,
\begin{align*}
R_{O}(K_{m},M)  &  <\frac{1}{2}MDC(K_{m},M)\\
R_{O}(K_{m},M)  &  =F_{g}(K_{m},M)\geq\varepsilon_{1}>\varepsilon^{\prime}\\
K  &  \in\overline{\mathcal{A}^{\infty}(k,\varepsilon^{\prime},D_{\sigma};M)}.
\end{align*}
\end{proof}

\begin{remark}
Different versions of the following lemma have been used by Whitney [W],
Cheeger and Gromov [CG], Gromov [Gr], Pugh [P] and others. Especially, the
proof of \ ''Hausdorff convergence implies Lipschitz convergence'' is along
the same lines. All versions known to the author are done in $\mathbf{R}^{N}$
to find a diffeomorphism between two Whitney embeddings. We include the
following version since it is in Riemannian manifolds with a uniform choice of
radius and an isotopy conclusion.
\end{remark}

\begin{lemma}
i. There exists $\rho(k,\varepsilon,D_{\varepsilon},M)>0$ such that for all
$K,L\in\mathcal{A}^{\infty}(k,\varepsilon,D;M)$ satisfying $L\subset
B(K,\rho;M)$ there exists a smooth isotopy between $K$ and $L$ in $B(K,\rho;M).$

ii. There exists $\rho^{\prime}(k,\varepsilon,D_{\varepsilon},M)>0$ such that
for all $K,L\in\mathcal{A}(k,\varepsilon,D;M)$ satisfying $L\subset
B(K,\rho^{\prime};M)$ there exists a continuous isotopy between $K$ and $L$ in
$B(K,\rho^{\prime};M)$ through $C^{1,1}$ embeddings of $L$.
\end{lemma}

\begin{proof}
i. Let $D^{\prime}=\overline{B}(D,\varepsilon)$ and $\varepsilon_{0}%
=\min(\varepsilon,\frac{1}{2}i(D^{\prime}))$. By Proposition 8(i), $\forall
K\in\mathcal{A}^{\infty}(k,\varepsilon,D;M),$ $\left\|  II^{K}\right\|  \leq
C_{0}$ and $\left|  Sect(K)\right|  \leq C_{1}.$ By Lemma 2, $\exists
d_{2}=\min(d_{1},\varepsilon_{0})>0$ such that and any $C^{2}$ curve
$\alpha:[0,d_{2}]\rightarrow D^{\prime}$ with $\left\|  \alpha^{\prime
}(s)\right\|  =1$ and $\left\|  \nabla_{\alpha^{\prime}}\alpha^{\prime
}\right\|  \leq C_{0}$ must satisfy $d_{M}(\gamma(s),\alpha(s))\leq\frac{s}%
{4},\forall s\in\lbrack0,d_{2}]$ where $\gamma(s)=\exp_{\alpha(0)}%
s\alpha^{\prime}(0)$. Given $p\in D^{\prime}$ and $v\in TM_{p}$, one can
naturally identify $T(TM_{p})_{v}\cong TM_{p}.$ The vector in $UT(TM_{p})_{v}$
corresponding to $u\in UTM_{p}$ under this identification will be denoted by
$u^{\prime},$ and let $u^{\prime\prime}=(d$ $\exp_{p})_{v}(u^{\prime}).$

Claim 1: $\exists d_{3}>0$ such that $d_{3}\leq d_{2}$ and $\forall p\in
D^{\prime},\forall u\in UTM_{p},\forall v\in TM_{p},\left\|  v\right\|  \leq
d_{3},$ one must have $d(\exp_{p}d_{2}u,\exp_{q}\frac{d_{2}u^{\prime\prime}%
}{\left\|  u^{\prime\prime}\right\|  })\leq\frac{d_{2}}{4}$ where $q=\exp
_{p}v,$ and $u^{\prime},u^{\prime\prime}$ are defined as above. Suppose that
such $d_{3}$ does not exist, then by using compactness, extract a subsequence
$p_{m}\rightarrow p_{0},$ $\left\|  v_{m}\right\|  \rightarrow0, $
$q_{m}\rightarrow q_{0},$ $u_{m}\rightarrow u_{0},$ $u_{m}^{\prime\prime
}\rightarrow u_{0}^{\prime\prime},$ with $d(\exp_{p_{m}}d_{2}u_{m},\exp
_{q_{m}}\frac{d_{2}u_{m}^{\prime\prime}}{\left\|  u_{m}^{\prime\prime
}\right\|  })>\frac{d_{2}}{4}$. By continuity and $d(\exp_{p})_{0}=Id,$ one
obtains $p_{0}=q_{0}$ and $u_{0}^{\prime\prime}=u_{0}.$ Then $d\left(
\exp_{p_{m}}d_{2}u_{m},\exp_{q_{m}}\frac{d_{2}u_{m}^{\prime\prime}}{\left\|
u_{m}^{\prime\prime}\right\|  }\right)  \rightarrow0$ which leads to a
contradiction. Thus, Claim 1 holds.

Set $\rho=\min\left(  \frac{d_{2}}{3},\frac{i_{0}}{4},d_{3}\right)  $ where
$0<i_{0}\leq i(K),\forall K\in\mathcal{A}^{\infty}(k,\varepsilon,D;M),$ by
Proposition 8(iv). Obviously, $\rho<\varepsilon_{0}<i(K,M).$

Let $K$ and $L$ be given as in the hypothesis. Define $E_{p}=\exp_{p}%
^{N}(B(0,\rho;NK_{p})),\forall p\in K,$ which are $C^{\infty}$ $(n-k)$
dimensional submanifolds of $M.$ $K^{k}$ is obviously transversal to all
$E_{p}.$

Claim 2. If $E_{p}\cap L\neq\emptyset,$ then $E_{p}$ intersects $L$
transversally at finitely many points. Suppose not: $\exists u^{\prime\prime
}\in T(E_{p})_{q}\cap TL_{q}-\{0\}$ for some $p\in K$ and $q\in E_{p}\cap L,$
since $\dim L+\dim E_{p}=\dim M.$ $\exists v\in NK_{p}\subset TM_{p}$ such
that $\left\|  v\right\|  <\rho\leq d_{3}$ and $\exp_{p}v=q.$ Let $u^{\prime
}=(d(\exp_{p})_{v})^{-1}(u^{\prime\prime})$ and adjust the length of
$u^{\prime\prime}$ so that $\left\|  u^{\prime}\right\|  =1.$ Find $u\in
UTM_{p}\cong UT(TM_{p})_{v}$ corresponding to $u^{\prime}.$ Then one has $u\in
UNK_{p}$, since $u^{\prime\prime}\in T(E_{p})_{q}$ and $E_{p}\subset\exp
_{p}(UNK_{p}).$ Set $a_{0}=\exp_{p}ud_{2}.$%
\begin{align*}
d_{2}  &  \leq\varepsilon_{0}<i(K,M)\\
B(a_{0},d_{2};M)\cap K  &  =\emptyset\\
d_{M}(a_{0},K)  &  =d_{2}\\
d_{M}\left(  a_{0},\exp_{q}\frac{d_{2}u^{\prime\prime}}{\left\|
u^{\prime\prime}\right\|  }\right)   &  \leq\frac{d_{2}}{4}\text{ by choice of
}d_{3}\geq\left\|  v\right\|  =d(p,q)\\
d_{M}\left(  \exp_{q}\frac{d_{2}u^{\prime\prime}}{\left\|  u^{\prime\prime
}\right\|  },\exp_{q}^{\mathbf{L}}\frac{d_{2}u^{\prime\prime}}{\left\|
u^{\prime\prime}\right\|  }\right)   &  \leq\frac{d_{2}}{4}\text{ by choice of
}d_{2}\text{, Lemma 2 and }\left\|  II^{L}\right\|  \leq C_{0}\\
d_{M}\left(  \exp_{q}^{\mathbf{L}}\frac{d_{2}u^{\prime\prime}}{\left\|
u^{\prime\prime}\right\|  },K\right)   &  \geq\frac{d_{2}}{2}>\rho
\end{align*}
The last assertion contradicts with $L\subset B(K,\rho;M).$ Hence, $\forall
q\in E_{p}\cap L,$ $T(E_{p})_{q}\cap TL_{q}=\{0\}$ and Claim 2 holds.

Since $K$ is smooth and $\rho<i(K,M),$ $\Psi=\left(  exp^{N}|B(0,\rho
,NK)\right)  $ is a diffeomorphism of $B(0,\rho,NK)$ onto $B(K,\rho,M).$
Define $\Pi:B(K,\rho,M)\rightarrow K$ by $\Pi^{-1}(p)=E_{p}.$ $\Pi$ is a
smooth submersion onto $K.$ By Claim 2, $\Pi|L:L\rightarrow K$ is a maximal
rank map. Since $L$ is compact, $K$ is connected, and $\dim L=\dim K,$ $\Pi|L$
must be onto and a covering map. Let $q_{1},q_{2}\in L$ be such that
$\Pi(q_{1})=\Pi(q_{2})=p.$%
\begin{align*}
p  &  \in B(q_{1},\rho,M)\cap B(q_{2},\rho,M)\\
q_{1}  &  \in L\cap B(q_{2},2\rho,M)\subset B(q_{2},3\rho,L)
\end{align*}
by $3\rho\leq\min(d_{2},i_{0})$, Lemma 2(i) and $\left\|  II^{L}\right\|  \leq
C_{0}.$ Let $\gamma$ be a normal minimal geodesic of $L$ from $q_{1}$ and
$q_{2}.$ The loop $\Pi(\gamma)$ is contractible in $K,$ since
\begin{align*}
d_{M}(\Pi(\gamma(t)),p)  &  \leq d_{M}(\Pi(\gamma(t)),\gamma(t))+d_{M}%
(\gamma(t),q_{j})+d_{M}(p,q_{j})\leq\frac{7}{2}\rho<i_{0}\leq i(K)\\
\text{where }j  &  =1\text{ for }t\leq\frac{3\rho}{2},\text{ and }j=2\text{
otherwise.}%
\end{align*}
By the homotopy lifting property, $q_{1}=q_{2}.$ Consequently, $\Pi|L$ is a
diffeomorphism of $L$ onto $K,$ and $\forall p\in K,E_{p}\cap L$ consists only
one point. Hence, $\Psi^{-1}(L)$ is a smooth section of the normal bundle
$B(0,\rho,NK)$ transverse to the fibers $NK_{p}$.The same is true for
$t_{0}\cdot\Psi^{-1}(L),$ $\forall t_{0}\in\lbrack0,1].$ Define $\Omega
:L\times\lbrack0,1]\rightarrow M$ by $\Omega(q,t)=\Psi(t\cdot\Psi^{-1}(q)).$
Obviously, $\Omega$ is a smooth map, $\Omega(q_{0},t)$ is the minimal geodesic
between $\Pi(q_{0})$ and $q_{0},$ and $\Omega(.,t_{0}) $ is a smooth embedding
of $L$ into $M,$ for all $t_{0}\in\lbrack0,1].$

ii. Choose $\varepsilon^{\prime}$ with $\mathcal{A}(k,\varepsilon
,D;M)\subset\overline{\mathcal{A}^{\infty}(k,\varepsilon^{\prime},D_{\sigma
};M)}$ and $\rho^{\prime}(k,\varepsilon,D_{\varepsilon},M)=\rho(k,\varepsilon
^{\prime},D_{\varepsilon},M).$ Consider any $K,L\in\mathcal{A}(k,\varepsilon
,D;M)$ satisfying $L\subset B(K,\rho^{\prime};M).$ By using Proposition 11,
find smooth approximations $K^{\prime}$ and $L^{\prime}$ with $d_{C^{1}%
}(K,K^{\prime})<\sigma$ and $d_{C^{1}}(L,L^{\prime})<\sigma$ for sufficiently
small $\sigma$ to secure $L^{\prime}\subset B(K^{\prime},\rho^{\prime};M).$
Recall that we constructed the smooth approximations $K^{\prime}$ by using
mollifiers locally in coordinate systems. One can construct the obvious
''vertical'' isotopies between the graphs: $(1-t)f(x)+th^{\delta}(x)$ for each
local smoothing and then push them forward into $M$ by the coordinate maps. By
applying these finitely many isotopies successively, one can construct an
isotopy between $K$ and $K^{\prime}.$ Similarly, one constructs an isotopy
between $L$ and $L^{\prime},$ and combining all one obtains an isotopy between
$K$ and $L.$ Other than the times of attachment of successive isotopies
constructed by using different local graphs or part (i), the isotopy is
$C^{1,1},$ and at any fixed time $t_{0}$ \ the embedding of $L$ is $C^{1,1}.$
\end{proof}

\subsection{Proof of Theorem 2}

\begin{proof}
We will take subsequences for several times, to simplify the notation all
subsequences will be denoted by the same index $m$, and $\forall m$ means
within the last chosen subsequence. The letter ''i'' appearing as a subindex
such as in $g_{ij}$ never means injectivity radius as in $i(K,M),$
$i(D^{\prime})$ or $i_{0}.$

i. By Proposition 11, $\exists\delta>0$ such that $\mathcal{D}(k,\varepsilon
,D;M)\subset\overline{\mathcal{A}^{\infty}(k,\delta,D_{\varepsilon};M)}$ in
$C^{1}$ topology. By Proposition 9, $\mathcal{A}^{\infty}(k,\delta
,D_{\varepsilon};M)$ has finitely many diffeomorphism types, and hence, the
same is true for $\mathcal{D}(k,\varepsilon,D;M)$. The finiteness of isotopy
classes will follow (ii) and Lemma 5(ii).

ii. Let a sequence $\{(K_{m},M)\}_{m=1}^{\infty}$ in $\mathcal{D}%
(k,\varepsilon,D;M)$ be given. By Proposition 11, $\forall m\in\mathbf{N}^{+}%
$, choose smooth submanifolds $\left(  L_{m},M\right)  \in\mathcal{A}^{\infty
}(k,\delta,D_{\varepsilon};M)$ such that $d_{C^{1}}(K_{m},L_{m})<\frac{1}{m}$.
Choose a subsequence so that all $L_{m}$ are diffeomorphic to a fixed
$C^{\infty}$ manifold $L,$ by the finiteness of diffeomorphism types. Hence,
$\forall m\in\mathbf{N}^{+},$ there are $C^{\infty}$ embeddings $e_{m}%
:L\rightarrow(M,g_{0})$ such that $L_{m}=e_{m}(L)$ and Riemannian metrics
$g_{m}=e_{m}^{\ast}g_{0}$ are $C^{\infty}$ on $L.$ By the intrinsic form of
Gromov's Compactness Theorem, as it was stated in [Pe, Thm. 4.4], there exists
a subsequence $g_{m}\rightarrow g_{\infty}$ in $C^{1}$ sense on $L$ with
respect to harmonic coordinates$,$ where $g_{\infty}$ is a $C^{1,\alpha}$
Riemannian metric on $L.$

We will show below that there exists an isometric embedding $e_{\infty
}:(L,g_{\infty})\rightarrow(M,g_{0})$ such that $e_{m}\rightarrow e_{\infty}$
in $C^{1}$ sense. Let $D^{\prime}=\overline{B}(D,\varepsilon)$ and $\delta
_{0}=\min(\delta,\frac{1}{2}i(D^{\prime}))$. Choose a finite collection of
points $p_{\alpha}$ such that $\{B(p_{\alpha},\delta_{0};M):\alpha
=1,...,\alpha_{0}\}$ covers $D^{\prime}$ and define $\varphi_{\alpha}:=\left(
\exp_{p_{\alpha}}^{M}|\overline{B(0,\delta_{0};TM_{p_{\alpha}})}\right)
^{-1}.$ There exists a Lebesgue number $c_{9}>0$ for the covering
$\{B(p_{\alpha},\delta_{0};M):\alpha=1,...,\alpha_{0}\},$ that is$:$ $\forall
q\in D^{\prime},\exists\alpha(q)$ such that $B(q,c_{9};D^{\prime})\subset
B(p_{\alpha(q)},\delta_{0};M).$ Let $r_{1}=\min(\frac{i_{0}}{4},\frac{c_{9}%
}{2})$ where $i(L,g_{m})\geq i_{0}>0,\forall m,$ by Proposition 8(iv).

By following [Pe, Thm. 4.4], \ for sufficiently large $m\in\mathbf{N}^{+}%
\cup\{\infty\},$ choose a finite open cover of $L$ by balls $\{B(q_{s}%
,r;(L,g_{m})):s=1,...,s_{0}(n,d_{0},v_{0},C_{1})\}$ such that the harmonic
coordinates of [JK] exist on $B(q_{s},2r;(L,g_{m}))$ for some $r=r(n,d_{0}%
,v_{0},C_{1}))\in(0,r_{1}].$ By [JH] and [Pe], the components of the metrics
in the harmonic coordinates satisfy
\[
\left(  g_{m}\right)  _{ij}\rightarrow\left(  g_{\infty}\right)  _{ij}\text{
in }C^{1}\text{ sense as }m\rightarrow\infty.\text{\qquad}(2.1)
\]
$\forall s=1,...,s_{0},$ the sequence $\{e_{m}(q_{s})\}_{m=1}^{\infty}$ is in
compact $D^{\prime}.$ By taking a subsequence, assume that $\forall
s=1,...,s_{0},$ $e_{m}(q_{s})\rightarrow z_{s}$, as $m\rightarrow\infty$ for
some $z_{s}\in D^{\prime}.$

Fix $s=1$. Set $\psi_{m}:B(q_{1},2r;(L,g_{m}))\rightarrow\mathbf{R}^{k}$ to be
the harmonic coordinates and $V=\overline{B(q_{s},r;(L,g_{\infty}))}.$ By the
construction of harmonic coordinates [JH], $(2.1)$, and $e_{m}$ being
isometric embeddings, there exists a compact $W\subset\mathbf{R}^{k}$ such
that for sufficiently large $m,$ all of the following holds.
\begin{align*}
V  &  \subset B(q_{1},2r;(L,g_{m}))\\
\psi_{m}\left(  V\right)   &  \subset int(W)\subset W\subset\psi_{m}\left(
B(q_{1},2r;(L,g_{m}))\right) \\
e_{m}(B(q_{1},2r;(L,g_{m})))  &  \subset B(z_{1},2r;D^{\prime})\subset
B(p_{\alpha(1)},\delta_{0};M)
\end{align*}
Define the following $C^{\infty}$ functions and vector fields for the last
subsequence:
\begin{align*}
h_{m}(y_{1},y_{2},..,y_{k})  &  =\varphi_{\alpha(1)}\circ e_{m}\circ\psi
_{m}^{-1}:W\subset\mathbf{R}^{k}\rightarrow TM_{p_{\alpha(1)}}\cong
\mathbf{R}^{n}\\
Y_{i}^{m}  &  =\left(  \psi_{m}^{-1}\right)  _{\ast}\left(  \frac{\partial
}{\partial y_{i}}\right)  \text{ and }Z_{i}^{m}=\left(  e_{m}\right)  _{\ast
}Y_{i}^{m}\
\end{align*}
All of the estimates below are uniformly on $W$, or the corresponding domains
by $\psi_{m}^{-1}$ and $e_{m}$, for all present indices $i,j,l,$ and $m$ of
the last chosen subsequence when limit is not taken.
\[
\left\langle Y_{i}^{m},Y_{j}^{m}\right\rangle _{g_{m}}=\left(  g_{m}\right)
_{ij}\rightarrow\left(  g_{\infty}\right)  _{ij}\text{ as }m\rightarrow
\infty\text{\qquad}(2.2)
\]
Since $g_{\infty}$ is a non-degenerate metric and $e_{m}$ are isometric
embeddings, there exists constants $c_{10}$, $c_{11}$and $c_{12}$ such that
\begin{align*}
0  &  <c_{10}\leq\left\|  Z_{j}^{m}\right\|  _{g_{0}}=\left\|  Y_{j}%
^{m}\right\|  _{g_{m}}\leq c_{11}<\infty\text{\qquad}(2.3)\\
Z_{l}^{m}\left\langle Z_{i}^{m},Z_{j}^{m}\right\rangle _{g_{0}}  &  =Y_{l}%
^{m}\left\langle Y_{i}^{m},Y_{j}^{m}\right\rangle _{g_{m}}=\frac
{\partial\left(  g_{m}\right)  _{ij}}{\partial y_{l}}\rightarrow\frac
{\partial\left(  g_{\infty}\right)  _{ij}}{\partial y_{l}}\text{ as
}m\rightarrow\infty\text{\qquad}(2.4)\\
\left\|  Z_{l}^{m}\left\langle Z_{i}^{m},Z_{j}^{m}\right\rangle _{g_{0}%
}\right\|   &  =\left\|  Y_{l}^{m}\left\langle Y_{i}^{m},Y_{j}^{m}%
\right\rangle _{g_{m}}\right\|  \leq c_{12}<\infty\text{\qquad}(2.5)\\
\left\|  \left\langle \nabla_{Y_{l}^{m}}^{m}Y_{i}^{m},Y_{j}^{m}\right\rangle
_{g_{m}}\right\|   &  =\frac{1}{2}\left\|  Y_{l}^{m}\left\langle Y_{i}%
^{m},Y_{j}^{m}\right\rangle _{g_{m}}-Y_{i}^{m}\left\langle Y_{l}^{m},Y_{j}%
^{m}\right\rangle _{g_{m}}+Y_{j}^{m}\left\langle Y_{i}^{m},Y_{l}%
^{m}\right\rangle _{g_{m}}\right\| \\
&  \leq\frac{3}{2}c_{12}\text{\qquad}(2.6)
\end{align*}
where $\nabla^{m}$ denotes the connection of $(K,g_{m}).$ Let $\nabla^{M}$
denote the connection of $(M,g_{0}).$ Since $e_{m}$ are isometric embeddings,
$L_{m}$ are $C^{\infty}$ submanifolds and $R_{O}(L_{m},M)\geq\delta,\exists
c_{13},c_{14}$ such that
\begin{align*}
\left\|  \left\langle \nabla_{Z_{l}^{m}}^{M}Z_{i}^{m},Z_{j}^{m}\right\rangle
_{g_{o}}\right\|   &  \leq\frac{3}{2}c_{12}\text{ \qquad}(2.6^{\prime})\\
\left\|  \left\langle \nabla_{Z_{j}^{m}}^{M}Z_{j}^{m},\overrightarrow
{n}\right\rangle _{g_{o}}\right\|   &  \leq C_{0}(\delta)\left\|  Z_{j}%
^{m}\right\|  _{g_{0}}^{2}\leq c_{13},\text{ }\forall\overrightarrow{n}\in
UNL_{m}\text{ \qquad}(2.7)\\
\left\|  \left\langle \nabla_{Z_{l}^{m}}^{M}Z_{j}^{m},\overrightarrow
{n}\right\rangle _{g_{o}}\right\|   &  \leq\frac{3}{2}c_{13},\forall
\overrightarrow{n}\in UNL_{m}\text{ \qquad}(2.8)\\
\left\|  \nabla_{Z_{l}^{m}}^{M}Z_{j}^{m}\right\|  _{g_{o}}  &  \leq
c_{14},\forall\overrightarrow{n}\in UNL_{m}\text{ \qquad}(2.9)
\end{align*}
There exists $c_{15},c_{16}$ depending on $M,$ $p_{\alpha(1)},$ and
$\exp_{p_{\alpha(1)}}^{M}$such that
\[
\forall v\in UTM|\overline{B(p_{\alpha(1)},\delta_{0};M)},\text{\ }%
0<c_{15}\leq\left\|  \left(  \varphi_{\alpha(1)}\right)  _{\ast}(v)\right\|
_{\mathbf{R}^{n}}\leq c_{16}<\infty.
\]%
\[
\left\|  \frac{\partial h_{m}}{\partial y_{j}}\right\|  _{\mathbf{R}^{n}%
}=\left\|  \left(  h_{m}\right)  _{\ast}\left(  \frac{\partial}{\partial
y_{j}}\right)  \right\|  _{\mathbf{R}^{n}}=\left\|  \left(  \varphi
_{\alpha(1)}\right)  _{\ast}(Z_{j}^{m})\right\|  _{\mathbf{R}^{n}}\leq
c_{16}\cdot c_{11}\text{\ \qquad}(2.10)
\]%
\[
\left(  \varphi_{\alpha(1)}\right)  _{\ast}\left(  \nabla_{Z_{l}^{m}}^{M}%
Z_{j}^{m}\right)  =\sum\limits_{\gamma}\left(  \sum_{\beta,\eta}\left(
\frac{\partial h_{m}}{\partial y_{l}}\right)  _{\beta}\left(  \frac{\partial
h_{m}}{\partial y_{j}}\right)  _{\eta}\Gamma_{\beta\eta}^{\gamma}+\left(
\frac{\partial^{2}h_{m}}{\partial y_{l}\partial y_{j}}\right)  _{\gamma
}\right)  \frac{\partial}{\partial x_{\gamma}}%
\]
in the local coordinates $\varphi_{\alpha(1)},$ where $()_{\beta}$ denotes the
$\beta$th component in $\mathbf{R}^{n}.$ Hence, $\exists c_{17}(k,n,c_{11}%
,c_{16},g_{0},\varphi_{\alpha(1)})<\infty$ such that
\[
\left\|  \frac{\partial^{2}h_{m}}{\partial y_{l}\partial y_{j}}\right\|
_{\mathbf{R}^{n}}\leq c_{17}.\text{\ \qquad}(2.11)
\]
Since $W$ is compact, $e_{m}(q_{1})\rightarrow z_{1}$ as $m\rightarrow\infty$,
$(2.10)$ and $(2.11),$ there exists a subsequence of $h_{m}$ converging in
$C^{1}$ topology to a $C^{1,1}$ function $h_{\infty}$ over $W, $ by
Arzela-Ascoli Theorem. Hence, there exists a subsequence of $e_{m}$ converging
to a $C^{1,1}$ function over $B(q_{1},r;(L,g_{\infty})).$ By applying this
process on all finitely many $B(q_{s},r;(L,g_{\infty})),$ there exists a
subsequence $e_{m}\rightarrow e_{\infty}\in C^{1,1}$ on $L$ in $C^{1}$ topology.

$e_{\infty}$ is an immersion since $h_{\infty}$ is non-singular: $\forall
v\in\mathbf{R}^{n}-\{0\},$ and sufficiently large $m,$
\begin{align*}
\left\|  \left(  h_{m}\right)  _{\ast}(v)\right\|  _{\mathbf{R}^{n}}  &
=\left\|  \left(  \varphi_{\alpha(1)}\right)  _{\ast}\circ\left(
e_{m}\right)  _{\ast}\circ\left(  \psi_{m}^{-1}\right)  _{\ast}(v)\right\|
_{\mathbf{R}^{n}}\\
&  \geq c_{15}\left\|  \left(  \psi_{m}^{-1}\right)  _{\ast}(v)\right\|
_{g_{m}}\geq\frac{c_{15}}{2}\left\|  \left(  \psi_{m}^{-1}\right)  _{\ast
}(v)\right\|  _{g_{\infty}}>0.
\end{align*}
Suppose that $e_{\infty}$ is not one-to-one, $e_{\infty}(a)=e_{\infty}(b)$ for
some $a,b\in L.$ Let $A=\frac{1}{5}\min(d_{2},d(a,b;(L,g_{\infty})))$ with
$d_{2}$ of Proposition 8(iii). For sufficiently large $m,$ $d(a,b;(L,g_{m}%
))\leq\frac{5}{4}d(a,b;(L,g_{\infty})).$ As in the proof of Proposition
8(iii),
\[
B(e_{m}(a),A;M)\cap B(e_{m}(b),A;M)=\emptyset.
\]
This contradicts with $e_{m}(a)\rightarrow e_{\infty}(a)$ and $e_{m}%
(b)\rightarrow e_{\infty}(b)=e_{\infty}(a).$ Hence, $e_{\infty}$ is one-to-one.

$g_{m}=e_{m}^{\ast}g_{0}\rightarrow e_{\infty}^{\ast}g_{0}$ in $C^{0}$
topology, since $e_{m}\rightarrow e_{\infty}$ in $C^{1}$ topology. However,
$g_{m}\rightarrow g_{\infty}$ in $C^{1}$ topology on $L$ in harmonic
coordinates. Consequently, $e_{\infty}^{\ast}g_{0}=g_{\infty},$ i.e.,
$e_{\infty}:(L,g_{\infty})\rightarrow(M,g_{0})$ is an isometric embedding.

Since we have chosen $d_{C^{1}}(K_{m},L_{m})<\frac{1}{m}$ in the beginning of
the proof, the initial sequence $\{(K_{m},M)\}_{m=1}^{\infty}$ has a $C^{1}%
$-convergent subsequence $\{(K_{m_{j}},M)\}_{j=1}^{\infty}$ whose limit is
$(e_{\infty}(L),M):=(K_{0},M).$

$(K_{0},M)$ belongs to $\mathcal{D}(k,\varepsilon,D;M)$ by $K_{m_{j}}\subset
D$ and Propositions 1 and 2:
\[
\underset{j\rightarrow\infty}{\varepsilon\leq\lim\sup\text{ \ }}i(K_{m_{j}%
},M)\leq i(K_{0},M).
\]

iii. Given For every $(K,M)\in\mathcal{D}(k,\varepsilon,D;M)$ with the given
embedding $e:K\rightarrow(M,g_{0}),$ find smooth approximations $\left(
L_{m},M\right)  $ of $(K,M)$ in $\mathcal{A}^{\infty}(k,\delta,D_{\varepsilon
};M)$ such that $d_{C^{1}}(K,L_{m})<\frac{1}{m}$, repeat (ii) to show that
$e(K)=e_{\infty}(L)$ by uniqueness of limits. Of course, $g_{\infty}$ is
$C^{1,\alpha}$ $(\alpha<1)$ by [Pe] in harmonic coordinates [JK] and it is a
limit of $C^{\infty}$ Riemannian metrics of bounded curvature and injectivity
radius with respect to Lipschitz distance, [Ni], [Pe], [GW]. $(K,g_{\infty})$
is a $C^{1,\alpha}$ Alexandrov space [Ni] with a well defined exponential map,
[D3], [Pu].
\end{proof}

\begin{corollary}
By Proposition 2 and Theorem 2, there exists thickest submanifolds in every
nonempty diffeomorphism or isotopy class $\mathcal{E}\subset\mathcal{D}%
(k,\varepsilon,D;M)$:
\[
\exists(K_{0},M)\in\mathcal{E}\text{ such that }\forall(K,M)\in\mathcal{E,}%
\text{ }i(K_{0},M)\geq i(K,M)
\]
The same conclusion holds for any nonempty closed subset $\mathcal{E}$ of
$\mathcal{D}(k,\varepsilon,D;M).$ For example, subsets with fixed volume or diameter.
\end{corollary}

\subsection{$C^{1}-$Compactness for $K$ with many components}

Define $\mathcal{D}^{\ast}(k,\varepsilon,D;M)=\{(K,M):K\in C^{1,1},$ $\dim
K=k,$ $K\subset D,$ and $i(K,M)\geq\varepsilon\}$ where $K$ is not necessarily connected.

\begin{corollary}
\ \ \ \ 

i. The number of components of $K$ in $\mathcal{D}^{\ast}(k,\varepsilon,D;M)$
are uniformly bounded.

ii. $\mathcal{D}^{\ast}(k,\varepsilon,D;M)$ is sequentially compact in $C^{1}%
$-topology, and it has finitely many isotopy and diffeomorphism types.

iii. There exists a thickest submanifold in each nonempty isotopy class of
$\mathcal{D}^{\ast}(k,\varepsilon,D;M).$
\end{corollary}

\begin{proof}
i. Let $\varepsilon_{0}=\min(\varepsilon,\frac{1}{2}i(D_{\varepsilon})).$ For
each component $K^{\alpha}$ of $K,$ choose $p_{\alpha}\in K^{\alpha}.$
$B(K^{\alpha},\varepsilon;M)\cap B(K^{\beta},\varepsilon;M)=\emptyset$ for
$\alpha\neq\beta,$ since $i(K,M)\geq\varepsilon.$ Hence, $B(p_{\alpha
},\varepsilon_{0};M)\cap B(p_{\beta},\varepsilon_{0};M)=\emptyset$. By Croke
[Cr, Prop. 14], $\exists v_{1}(n,\varepsilon_{0})>0$ such that $v_{1}\leq
vol_{n}(B(p_{\alpha},\varepsilon_{0};M)),\forall\alpha$. Hence, $K$ has at
most $vol_{n}(D_{\varepsilon})/v_{1}$ components.

ii. Given any sequence $\{(K_{j},M)\}_{j=1}^{\infty}$ in $\mathcal{D}^{\ast
}(k,\varepsilon,D;M),$ choose a subsequence (by using same index $j$) where
the number of components of $K_{j}$ is constant, and enumerate the components
$K_{j}^{\alpha}.$ $\forall\alpha,$ $\varepsilon\leq i(K_{j},M)\leq
i(K_{j}^{\alpha},M).$ By Theorem 2, choose a subsequence where $K_{j}^{1}$
converges in $C^{1}$ topology, and choose its subsequence where $K_{j}^{2}$
converges in $C^{1}$ topology and so on. Hence, $\mathcal{D}^{\ast
}(k,\varepsilon,D;M)$ is sequentially compact in $C^{1}$-topology and a
subsequence $(K_{j},M)\rightarrow(K_{\infty},M)\in\mathcal{D}^{\ast
}(k,\varepsilon,D;M)$, by Proposition 2. $\exists j_{0}\forall j\geq
j_{0}\forall\alpha,$ $K_{j}^{\alpha}\subset B(K_{\infty}^{\alpha},\rho
^{\prime};M)$, for $\rho^{\prime}$ of Lemma 5 and hence $K_{j}^{\alpha}$ is
isotopic to $K_{\infty}^{\alpha}$ in $B(K_{\infty}^{\alpha},\varepsilon;M).$
These isotopies can be combined to give an isotopy of $K_{j}$ to $K_{\infty}$
without any self intersections $\forall j\geq j_{0},$ since $B(K_{\infty
}^{\alpha},\varepsilon;M)$ are mutually disjoint.

iii. This is an immediate consequence of (ii) and Proposition 2.
\end{proof}

\subsection{Normal curvatures and Thickness Formula in Euclidean Spaces}

$\exp_{p}:T(K,g_{\infty})_{p}\rightarrow(K,g_{\infty})$ is of class $C^{0,1},$
see [D3] and [Pu]. Even though the geodesics $\exp_{p}sv$ of $(K,g_{\infty})$
are $C^{2}$ in limit harmonic coordinates, [D3, 5.10.3], the corresponding
geodesics $e_{\infty}(\exp_{p}sv)$ in $(K,M)$ are $C^{1,1}.$ Hence
$\nabla_{\gamma^{\prime}}^{M}\gamma^{\prime}$ is defined almost everywhere in $s.$

\begin{definition}
We define the supremum of the ''absolute normal curvatures'' $\sup\kappa
_{N}(K)$ for a $C^{1,1}$ submanifold $K$ to be

$\sup\left\{  \left\|  \nabla_{\gamma^{\prime}}^{M}\gamma^{\prime}(s)\right\|
:\gamma:\mathbf{R}\rightarrow K\text{ is a geodesic of }K\text{ with }\left\|
\gamma^{\prime}\right\|  =1\text{ and }\nabla_{\gamma^{\prime}}^{M}%
\gamma^{\prime}(s)\text{ exists}\right\}  .$
\end{definition}

\begin{proposition}
For a $C^{1,1}$ submanifold $K^{k}$ of $\mathbf{R}^{n},$ $F_{g}(K,\mathbf{R}%
^{n})=\frac{1}{\sup\kappa_{N}(K)}.$ Hence,
\[
i(K,M)=\min\{\frac{1}{\sup\kappa_{N}(K)},\frac{1}{2}MDC(K)\}.
\]
\end{proposition}

\begin{proof}
The following is a basic result in $\mathbf{R}^{n},$ we refer to [D6,
Proposition 2] for an elementary proof.

Let $\gamma:I=(-\frac{\pi}{2\kappa},\frac{\pi}{2\kappa})\rightarrow
\mathbf{R}^{n}$ be with $\left\|  \gamma^{\prime}\right\|  \equiv1,$ $\left\|
\gamma^{\prime\prime}\right\|  \leq\kappa\neq0$ a.e. Then,

\qquad i. $\gamma\cap O_{\gamma(0)}(\gamma^{\prime}(0),\frac{1}{\kappa
};\mathbf{R}^{n})=\emptyset,$ and

\qquad ii. if $\gamma^{\prime\prime}(0)$ exists and $\left\|  \gamma
^{\prime\prime}(0)\right\|  =\kappa,$ then $\forall R>\frac{1}{\kappa}%
,\exists\delta>0$ such that $\gamma((0,\delta))\subset B(\gamma(0)+R\frac
{\gamma^{\prime\prime}(0)}{\left\|  \gamma^{\prime\prime}(0)\right\|  },R).$

Let $\kappa_{0}=\sup\kappa_{N}(K)$ and $\varepsilon=\frac{\pi}{2\kappa_{0}}.$
By (i), for any $p\in K,$ $v\in UTK_{p},$
\begin{align*}
\exp_{p}^{K}((-\varepsilon,\varepsilon)v)\cap O_{p}(v,\frac{1}{\kappa_{0}%
},\mathbf{R}^{n})  &  =\emptyset\\
B(p,\varepsilon;K)\cap O_{p}(\frac{1}{\kappa_{0}},K;\mathbf{R}^{n})  &
=\emptyset\\
F_{g}(p,K;\mathbf{R}^{n})  &  \geq\frac{1}{\kappa_{0}}\\
F_{g}(K;\mathbf{R}^{n})  &  \geq\frac{1}{\kappa_{0}}=\inf_{\gamma\text{
geodesic}}\frac{1}{\left\|  \gamma^{\prime\prime}\right\|  }%
\end{align*}
Suppose that $F_{g}(K;\mathbf{R}^{n})>\underset{\gamma}{\inf}\frac{1}{\left\|
\gamma^{\prime\prime}\right\|  }.$ Then, there exists a geodesic $\gamma$ and
$R$ such that $\gamma^{\prime\prime}(0)$ exists and $F_{g}(K;\mathbf{R}%
^{n})>R>\frac{1}{\left\|  \gamma^{\prime\prime}(0)\right\|  }.$ Then by (ii) above,

$\gamma((0,\delta))\subset B(\gamma(0)+R\frac{\gamma^{\prime\prime}%
(0)}{\left\|  \gamma^{\prime\prime}(0)\right\|  },R)$ which implies that
$R\geq F_{g}(\gamma(0),K;\mathbf{R}^{n})\geq F_{g}(K;\mathbf{R}^{n})$ by the
definition of $F_{g}.$ Hence, one obtains a contradiction. Consequently,
$F_{g}(K,\mathbf{R}^{n})=\frac{1}{\sup\kappa_{N}(K)}.$ The rest follows
Theorem 1.
\end{proof}

\section{Estimates on the number of Isotopy and Diffeomorphism types}

The number $\#(k,\varepsilon,D;M)$ of the different diffeomorphism classes and
isotopy classes of $C^{1,1}$ manifolds of $\mathcal{D}(k,\varepsilon,D;M)$ is
bounded above by a constructible constant in terms $n,k,\delta(\varepsilon)$
and $D_{\varepsilon}$ where $\mathcal{D}(k,\varepsilon,D;M)\subset
\overline{\mathcal{A}^{\infty}(k,\delta,D_{\varepsilon};M)}$ in $C^{1}$
topology. It is clear from the proofs of Propositions 8 and 11, and Lemmas 2
and 5 that $\rho=\rho(n,k,\delta(\varepsilon),C_{0},\left|  Sect(M)\right|
,i(D_{\varepsilon})).$ The dependence of $\delta$ on $\varepsilon$ relies on a
finite number of fixed coordinate charts of $M,$ in fact on their derivatives
up to second order. Using normal coordinates may bring in $\left\|  \nabla
R\right\|  ,$ but in harmonic coordinates, one can control them with only
$\left|  Sect(M)\right|  $ and $i(D_{\varepsilon}). $

The number of different diffeomorphism classes and isotopy classes in
$\mathcal{D}(k,\varepsilon,D;M)$ can be bounded in terms of $\rho
(\delta(\varepsilon))$ as follows. Take a minimal cover $\mathcal{B}%
=\{B(p_{\alpha},\rho/2):\alpha=1,...\Lambda_{0}\}$ of $D_{\varepsilon}$ by
open discs: $B(p_{\alpha},\rho/4)\cap B(p_{\beta},\rho/4)=\emptyset$ if
$\alpha\neq\beta.$ Then, $\Lambda_{0}\leq vol(D_{\varepsilon})/\min
(vol(B(p_{\alpha},\rho/4)))\leq c(n)vol(D_{\varepsilon})\rho^{-n}$ by volume
estimates of [Cr]. Define $\Phi(K)=\{\alpha:K\cap B(p_{\alpha},\rho
/2)\neq\emptyset\}.$ Any $K,L\in\mathcal{D}(k,\varepsilon,D;M)$ with
$\Phi(K)=\Phi(L)$ must satisfy $L\subset B(K,\rho;M),$ and hence, $K$ and $L$
are isotopic and diffeomorphic. Consequently, there are at most
$2^{c(n)vol(D_{\varepsilon})\rho^{-n}}$ distinct diffeomorphism classes and
isotopy classes of $C^{1,1}$ manifolds in $\mathcal{D}(k,\varepsilon,D;M).$

We calculate these estimates in $\mathbf{R}^{n}$ below. Let $D=\overline
{B(0,r,\mathbf{R}^{n})}$ and $i(K,\mathbf{R}^{n})\geq\varepsilon.$ Rescale the
metric so that $R=\frac{r}{\varepsilon},D=\overline{B(0,R,\mathbf{R}^{n})}$
and $i(K,\mathbf{R}^{n})\geq1.$ Let $\alpha_{n}=vol(S^{n}(1)).$
\begin{align*}
i(\mathbf{R}^{n})  &  =\infty,\text{ }Sect(\mathbf{R}^{n})=0,\text{ and }%
C_{0}=C_{1}=1\\
1  &  =\varepsilon\geq\min(\varepsilon^{\prime},\delta,\delta_{0})\approx1\\
d_{1}  &  =d_{2}=\frac{1}{2}\text{ and }d_{3}=\frac{1}{8}\\
d_{0}  &  \leq\frac{1}{2}(8R)^{n}\\
v_{0}  &  \geq\frac{(n-1)\alpha_{n}}{n\alpha_{n-k-1}}\cdot e^{1-n}%
\end{align*}%
\begin{align*}
i_{0}  &  \geq\min(\pi,\frac{\pi v_{0}}{\alpha_{n}}\sinh^{1-n}d_{0})\text{ by
[HK]}\\
&  \geq e^{-\frac{n}{2}(8R)^{n}},\text{ for }k\geq2,\\
i_{0}  &  \geq\pi,\text{ for }k=1
\end{align*}%
\begin{align*}
\rho &  =\min\left(  \frac{d_{2}}{3},\frac{i_{0}}{4},d_{3}\right)
=\frac{i_{0}}{4}\text{ for }k\geq2\\
\rho &  =\min\left(  \frac{d_{2}}{3},\frac{i_{0}}{4},d_{3}\right)
=d_{3}=\frac{1}{8},\text{ for }k=1\\
\Lambda_{0}  &  \leq\left(  \frac{4R}{\rho}\right)  ^{n}=\left(  \frac
{16R}{i_{0}}\right)  ^{n}\\
\#(k,\varepsilon,D;M)  &  \leq2^{\Lambda_{0}}%
\end{align*}

Almost all of the estimates are reasonable, except for $i_{0}$ and $v_{0}$ for
$k\geq2.$

\section{References}

[B] M. Berger, \textit{Une borne inferieure pour le volume d}$^{\prime}%
$\textit{une variete Riemannienne en fonction du rayon d}$^{\prime}%
$\textit{injectivite,} Ann. Inst. Fourier, Grenoble, \textbf{30} (1980), 259-265.

[Be] A. L. Besse, \textit{Manifolds all of whose geodesics are closed,
}Ergebnisse Math. Grenzgebiete, vol 93, Springer, Berlin 1978.

[BS] G. Buck and J. Simon, \textit{Energy and lengths of knots, }Lectures at
Knots 96, 219-234.

[Ch] J. Cheeger, \textit{Finiteness theorems for Riemannian manifolds,
}American Jour. of Math. \textbf{92 }(1970) 61-74.

[CE] J. Cheeger and D. G. Ebin, \textit{Comparison theorems in Riemannian
geometry, Vol 9, }North-Holland, Amsterdam, 1975.

[CG] J. Cheeger and M. Gromov, \textit{On characteristic numbers of complete
manifolds of bounded curvature and finite volume, }in Differential Geometry
and Complex Analysis, editors I. Chavel and H. M. Farkas, Springer-Verlag,
Berlin, 1985.

[Cn]\ S. S. Chern, \textit{Curves and surfaces in Euclidean space, }in Studies
in Global Geometry and Analysis, Math. Assoc. Amer., Prentice Hall, 1967, 16-56.

[Cr]\ C. Croke, \textit{Some isoperimetric inequalities and eigenvalue
estimates}, Ann. Sci. Ecole Norm. Sup., \textbf{13} (1980) 419-435.

[DoC]\ M. P. DoCarmo, \textit{Riemannian geometry, }Birkhauser, Cambridge,
Massachusetts, 1992.

[D1]\ O. C. Durumeric, \textit{Manifolds with almost equal diameter and
injectivity radius}, Jour. of Diff. Geom.\textit{, }\textbf{19} (1984) 453-474.

[D2]\ O. C. Durumeric, \textit{A generalization of Berger's almost 1/4 pinched
manifolds theorem I, }Bull. Am. Math. Soc., \textbf{12} (1985) 260-264.

[D3]\ O. C. Durumeric, \textit{A generalization of Berger's theorem on almost
1/4 pinched manifolds II, }Jour. Diff. Geom., \textbf{26} (1987) 101-139.

[D4]\ O. C. Durumeric, \textit{Finiteness theorems, average volume and
curvature}, Am. Jour. of Math., \textbf{111} (1989) 973-990.

[D5]\ O. C. Durumeric, \textit{Growth of fundamental groups and isoembolic
volume and diameter, }Proc. Am. Math. Soc.,\textbf{130 }(2002), no. 2, 585-590.

[D6]\ O. C. Durumeric, \textit{Local structure of ideal knots, }preprint.

[GM]\ O. Gonzales and H. Maddocks, \textit{Global curvature, thickness and the
ideal shapes of knots}, Proceedings of National Academy of Sciences,
\textbf{96 }(1999) 4769-4773.

[GWY]\ K. Grove, P. Petersen and J.-Y. Wu, \textit{Geometric finiteness
theorems via controlled geometry}, Invent. Math., \textbf{99} (1990), 205-213;
Erratum, \textbf{104} (1991) 221-222.

[GW]\ R. E. Greene and H. Wu, \textit{Lipschitz convergence of Riemannian
manifolds}, Pacific Jour. of Math., \textbf{131} (1988), 119-141

[GKM]\ D. Gromoll, W. Klingenberg and W. Meyer, \textit{Riemannsche Geometrie
im Grossen}, Lecture notes in Mathematics No. 55, Springer-Verlag, 2nd
Edition, 1975.

[Gr]\ M. Gromov redige par J. LaFontaine and P. Pansu, \textit{Structures
metriques pour les} \textit{varietes Riemaniennes}, Cedic/Ferdnan Nathan, 1981.

[HK]\ E. Heintze and H. Karcher, \textit{A general comparison theorem with
applications to volume estimates for submanifolds}, Ann. scient. Ec. Norm.
Sup. 4e serie, \textbf{11} (1978),451-470.

[JK]\ J. Jost and H. Karcher, \textit{Geometrische Methoden zur Gewinnung von
a-priori-schranken fur Harmonische Abbildingen}, Manuscripta Math.,
\textbf{40} (1982), 27-77.

[Ka]\ V. Katrich, J. Bendar, D. Michoud, R.G. Scharein, J. Dubochet and A.
Stasiak, \textit{Geometry and physics of knots}, Nature, \textbf{384} (1996) 142-145.

[K]\ A. Katsuda, \textit{Gromov's convergence theorem and its application},
Nagoya Math. Jour., \textbf{100} (1985), 11-48.

[L]\ Litherland, Unbearable thickness of knots, preprint.

[LSDR]\ A. Litherland, J Simon, O. Durumeric and E. Rawdon, \textit{Thickness
of knots}, Topology and its Applications, \textbf{91}(1999) 233-244.

[M]\ J. Milnor, \textit{Morse theory, }Annals of Math. Studies vol. 51,
Princeton Univ. Press, Princeton, 1973.

[N]\ A. Nabutovsky, \textit{Non-recursive functions, knots ''with thick
ropes'' and self-clenching ''thick'' hyperspheres}, Communications on Pure and
Applied Mathematics\textit{, }\textbf{48} (1995) 381-428.

[Ni]\ I. G. Nikolaev, \textit{Bounded curvature closure of the set of compact
Riemannian manifolds}, Bulletin of A.M.S.\textit{,}\textbf{24} (1) 171-176.

[Pe]\ S. Peters, \textit{Convergence of Riemannian manifolds}, Compositio
Math., \textbf{62} (1987), 3-16.

[Pu]\ C. C. Pugh, \textit{C}$^{1,1}$ \textit{conclusions in Gromov's theory},
Ergodic Theory Dynam. Systems, \textbf{7} (1987), no. 1, 133-147.

[Wa]\ F. W. Warner, \textit{Conjugate loci of constant order}, Annals of
Math., \textbf{86} (1967), 192-212.

[W]\ H. Whitney, \textit{Geometric integration theory}, Princeton Univ. Press,
Princeton, 1957.

[Y]\ T. Yamaguchi, \textit{Homotopy type finiteness theorems for certain
precompact families of Riemannian manifolds}, Proc. Amer. Math. Soc.,
\textbf{102}, No. 3 (1988) 660-666.
\end{document}